\documentclass[11pt,a4paper,american,reqno]{amsart}
\usepackage{amssymb}
\usepackage{amsmath}
\usepackage{graphicx}
\usepackage{enumerate}
\usepackage{physics}
\usepackage{epsf}
\usepackage{amsbsy}
\usepackage{amsfonts}
\usepackage{eucal}
\usepackage{amsthm}
\usepackage{secdot}
\usepackage{esint}
\usepackage{bbm}
\usepackage{varwidth}
\usepackage{tasks}
\usepackage{cite}
\usepackage[T1]{fontenc}
\usepackage[utf8]{inputenc}
\usepackage[colorlinks,pdfpagelabels,pdfstartview = FitH,bookmarksopen = true,bookmarksnumbered = true,linkcolor = red,plainpages = false,hypertexnames = false,citecolor = blue,pagebackref=false]{hyperref}
\usepackage{cleveref}

\theoremstyle{plain}
\newtheorem{thm}{Theorem}[section]
\newtheorem{lem}[thm]{Lemma}
\newtheorem{cor}[thm]{Corollary}
\newtheorem{prop}[thm]{Proposition}
\theoremstyle{definition}
\newtheorem{df}[thm]{Definition}
\newtheorem{assumption}[thm]{Assumptions}
\theoremstyle{remark}
\newtheorem{rk}[thm]{Remark}
\newtheorem{exmp}[thm]{Example}
\usepackage{amsmath, amssymb}
\usepackage[most]{tcolorbox}
\renewcommand{\d}{\:\! \mathrm{d}}

\allowdisplaybreaks
\newcommand{\Lp}{L^p(\Omega)}
\newcommand{\Wkp}{W^{1,p}_0(\Omega)}
\newcommand{\mrp}{M^r(p;\Omega)}

\usepackage{xcolor}
\long\def\symbolfootnote[#1]#2{\begingroup
\def\thefootnote{\fnsymbol{footnote}}\footnote[#1]{#2}\endgroup}
\numberwithin{equation}{section}

\begin{document}
\title[Properties of principal eigenvalue] {On the Principal Eigenvalue of the $(p,A)$-Laplacian with Sign-Changing Morrey Potentials: Complete Continuity, Optimization, and Differentiability}
\author{Pradeep Gupta, Ratan Kr. Giri and Sekhar Ghosh* }
\address[Pradeep Gupta]{Department of Mathematics, The LNM Institute of Information Technology, Jaipur, India - 302031
}
\email{pg48350@gmail.com / 24PMT003@lnmiit.ac.in}
\address[Ratan Kr. Giri]{Department of Mathematics, The LNM Institute of Information Technology, Jaipur, India - 302031
}
\email{: giri90ratan@gmail.com / ratan.giri@lnmiit.ac.in }
\address[Sekhar Ghosh]{Department of Mathematics, National Institute of Technology Calicut, Kozhikode, Kerala, India - 673601}
\email{sekharghosh1234@gmail.com / sekharghosh@nitc.ac.in}
% \thanks{{\em 2020 Mathematics Subject Classification:} 35J60; 35J70; 35J92; 35P30}
\thanks{*Corresponding author: S. Ghosh}
\begin{abstract}
We establish a unified framework to study the continuity and complete continuity of the principal eigenvalue functional $\lambda_1(V,m)$ associated with the following $(p,A)$-Laplacian operator involving both sign-changing potential $V$, and positive weight $m$ in a Morrey space: 
\begin{align*}
		-\Delta_{p,A}(u) + V|u|^{p-2}u &= \lambda m|u|^{p-2}u \quad \text{in}\,\, \Omega, \\
		u 	&= 0  \quad\quad\quad \quad\quad \text{on}\,\, \partial\Omega, 
	\end{align*}
where $1<p<N$, $\lambda\in\mathbb{R}$ is a parameter, and  $\Omega \subset \mathbb{R}^N (N \geq 2)$ is a bounded domain with boundary $\partial \Omega$. As a consequence, we investigate the existence and uniqueness of optimal principal eigenvalues $\lambda_1(V,m)$ over a weakly compact subset of the Morrey space. Finally, we derive a formula for the derivative of $\lambda_1(V,m)$ with respect to regular perturbations of potential $V$, and positive weight $m$ within a closed subspace of Morrey space. 

Our results develop analytic tools to overcome the lack of separability and reflexivity in Morrey spaces to investigate this eigenvalue problem. In particular, the study is new even for the linear case $p=2$.  We also emphasize that the tools established are applicable for a wide class of PDEs with coefficients in Morrey spaces.  
%Owing to the strict inclusion of Lebesgue space $L^p$ in Morrey spaces, these results generalize previous findings in the literature to Morrey spaces that lack both reflexivity and separability.
\end{abstract}
\maketitle
\section{Introduction}
The study of the eigenvalue problems for nonlinear elliptic operators involving a potential $V$, specifically for the operator $-\Delta_{p}(u)+V|u|^{p-2}u$ constitutes a prominent area of research in partial differential equations (PDEs), due to their numerous applications in various fields such as quantum mechanics, scattering theory, stability of bulk matter, image processing, and related areas; for example, see \cite{simon2007,simon2022,gilboa18,bonnaillie2015}. The presence of the potential $V$ significantly influences the operator's coercivity, compactness, and the spectral properties. Over the past few decades, several authors have investigated the existence,  qualitative properties, and optimization of eigenvalues for the operator $-\Delta_{p}(u)+V|u|^{p-2}u$ under diverse conditions on the potential $V$, for instance, the reader can refer to \cite{Pinchover16,Pezzo06,Giri23} and the references therein. Focusing on the linear case $p = 2$ without potential ($V = 0$), several works \cite{Noris2026,Henrot2006,Henrot2019,Veronica2025} have also discussed spectral shape optimization for the first eigenvalue of both the Dirichlet and Neumann Laplacian problems in the presence of an obstacle. These studies have been further extended to weighted eigenvalue problems, incorporating various assumptions on potential functions as well as non-constant weights \cite{Binding95,Cuesta01,Cuesta09,Otani88,anoop09}. 

In this context, we investigate several fundamental properties of the functional associated with the principal eigenvalue $\lambda_1(V, m)$, alongside its applications to optimization problem for the $(p,A)$-Laplacian operator involving sign-changing potential $V$, and positive weight $m$ of the form:
\begin{equation}
	\left\{
	\begin{aligned}
		-\Delta_{p,A}(u) + V|u|^{p-2}u &= \lambda m|u|^{p-2}u \quad \text{in}\,\, \Omega, \\
		u 	&= 0  \quad  \text{on}\,\, \partial\Omega, \end{aligned} \label{EVP}
	\right. 
\end{equation}
where $1<p<N$, $\lambda \in \mathbb{R}$, and  $\Omega \subset \mathbb{R}^N (N \geq 2)$ is a bounded domain with boundary $\partial \Omega$. Here,  $\Delta_{p,A}(u) := \operatorname{div}(|\nabla u|^{p-2}_A A\nabla u)$ denotes the $(p,A)$-Laplacian operator, where the matrix $A$ satisfies certain conditions. The potential $V$, and weight function $m$ belong to Morrey spaces satisfying certain assumptions. Throughout this paper, the matrix $A$, potential $V$, and weight function $m$ satisfy the following conditions:
\begin{assumption}\label{assump}  
		\begin{itemize}
			\item[(H1)] $A=(a_{ij}(x))_{i,j=1}^N \in L^\infty(\Omega;\mathbb{R}^{N\times N})$ is a symmetric matrix satisfying the uniform ellipticity condition over $\Omega$ i.e., there exists  $\Theta_\Omega > 0$ such that
            $$ \Theta_{\Omega} |\eta|\le |\eta|_{A} \le \Theta_{\Omega}^{-1} |\eta|, \quad \text{for all }\, \eta =(\eta_1, \eta_2, \cdots, \eta_N) \in \mathbb{R}^N \text{ and a.e. } x\in \Omega,$$
			where $\displaystyle{|\eta|^2_{A}:= \sum_{i,j=1}^N a_{ij}(x)\eta_i\eta_j}$. 
			\item[(H2)] For $r> \frac{N}{p}$,  $(V,m) \in M^r(p;\Omega) \times \mrp$   with $m>0 $ a.e. on $\Omega$, where $\mrp$ is denotes the special class of  Morrey spaces (see Definition \ref{splorrey}).
		\end{itemize}
	\end{assumption}
    The existence and optimization theory for the principal eigenvalue associated with the problem of the type \eqref{EVP} has been discussed by several researchers. In their pioneering work, Harrell \cite{Harrell1984} and later, Ashbaugh and Harrell \cite{Ashbaugh1987} established the existence and uniqueness of both the maximizers and minimizers for the principal eigenvalue of $-\Delta + V$ over a closed, convex, and bounded subset of $L^p(\Omega)$, assuming that potential $V \in L^p(\Omega)$ with $p\geq 1$. These findings were further extended to the nonlinear case involving the $p-$Laplacian operator $\Delta_p$ with $V \in L^p(\Omega)$ by Pezzo and Bonder \cite{Pezzo06}. We also refer to the work of Marras et al. \cite{Marras13}, where the authors studied the optimization problem for the principal eigenvalue of \eqref{EVP} for the case when $A=I_N$, $V =0$, and weight $m=\chi_D$, where $D \subset \Omega$ is a variable subset with fixed measure $\alpha$.

In addition to the analysis of the optimization problem for the principal eigenvalue $\lambda_1(V,m)$ over a closed, convex, and bounded subset, several studies have also investigated this problem over the rearrangement class of a given pair $(V_0, m_0)$ lying in suitable Lebesgue spaces. An interesting contribution in this direction was done by Cuccu et al. \cite{Cuccu08}, who studied the case $A=I_N$, $V=0$, and assuming $m$ lies in the rearrangement class of a fixed positive weight $m_0 \in L^\infty(\Omega)$. Later,  Pezzo and Bonder \cite{Pezzo10} generalized these results to the case $A=I_N$, with $(V,m)$ belonging to the rearrangement class of  a fixed pair $(V_0,m_0) \in L^r(\Omega) \times L^r(\Omega)$, where $r > \frac{N}{p}$ if $ 1 < p \leq N$ and $ r = 1$ if $ p > N$. Further, Biswas et al. \cite{Biswas23} extended this analysis to the pair $(V_0,m_0)\in L^r(\Omega) \times L^r(\Omega)$ with critical exponent $r = \frac{N}{p}$ when  $1<p < N$, and $r \in (1,\infty)$ when $p \geq N$. 

\subsection{The main contribution:} Following our recent work  Gupta et al. \cite{Gupta2026}, where we established the existence and qualitative properties of the principal eigenvalue $\lambda_1(V,m)$ for the problem \eqref{EVP} under Assumption \ref{assump}, the present paper further investigates several properties of the associated functional $\lambda_1(V,m)$ and its optimization over a suitable subset of $\mrp \times \mrp$. Since the Lebesgue space $L^r(\Omega)$ is strictly contained in the Morrey space $\mrp$ for $r > \frac{N}{p}$ with $1 < p < N$ (see Remark \ref{Remark} and Definition \ref{splorrey}), the main goal of this paper is to investigate the effect of Morrey potentials to the problem \eqref{EVP}. Specifically, our work extends and advances the results of \cite{Cuccu08,Pezzo06,Pezzo10,Biswas23} to a broader setting, i.e.,  we assume the pair $(V, m)\in \mrp \times \mrp$  with $m >0$ a.e. in $\Omega$. Unlike the Lebesgue spaces $L^r(\Omega)$ for $r > 1$, working within the Morrey space $\mrp$ introduces significant technical challenges. Specifically, $\mrp$ is neither reflexive nor separable, and smooth functions with compact support are not dense in it, i.e., $\overline{C_c^\infty(\Omega)} \neq \mrp$. To overcome these challenges, we develop here unified tools to prove the continuity and complete continuity of the principal eigenvalue functional $\lambda_1(V,m)$ associated with problem \eqref{EVP}. As a consequence, we study the following optimization problems:
$$T_1:=\inf\{ \lambda_1(V,m) : (V,m) \in \mathcal{W}_C\},   \text{ and } T_2 := \sup\{ \lambda_1(V,m) : (V,m) \in \mathcal{W}_C\},$$ where $\mathcal{W}_C$ is a weakly compact subset of $\{ (V,m) \in \mrp \times \mrp : m > 0 \text{ a.e. in } \Omega \}$, and establish the existence of a minimizer for $T_1$ and a maximizer for $T_2$ in $\mathcal{W}_C$, respectively. In addition, when $\mathcal{W}_C$ is strongly convex, we prove that the maximizer of $\lambda_1(V,m)$ is unique, whereas the minimizer may fail to be unique; a fact we illustrate by constructing a counterexample. Moreover, we derive the derivative formulae for $\lambda_1(V, m)$ with respect to regular perturbations of potential $V$ and positive weight $m$ within a closed subspace  $M^r_0(p;\Omega):= \overline{C_c^\infty(\Omega)}$ of the Morrey space $\mrp$. It is worth noting that our established results remain valid for the range $p>N$. In this range, the special Morrey space $\mrp$ coincides with $L^1(\Omega)$ (see Definition 2.3 of $\mrp$, and the crucial Morrey-Adams Theorem 2.4 in Pinchover and Psaradakis \cite{Pinchover16}), and the proofs are analogous. However, the critical range $p=N$ requires a separate treatment. %and will be addressed in a forthcoming work.

In respect of the approximation of eigenvalues via the Finite Element Method (see \cite{bonder2007, henrot2005} and the relevant references therein),  we expect that these properties and formulas for $\lambda_1(V,m)$ can serve as a suitable framework for the numerical computation of optimal eigenvalues. 

\medskip

 This paper is arranged as follows. Section $2$ is devoted to reviewing essential results about Morrey spaces, as well as the principal eigenvalue $\lambda_1(V,m)$ of problem \eqref{EVP} and its qualitative properties. In Section $3$, we establish several crucial lemmas concerning the convergence of integrals of the form $\int_\Omega V_n |u_n|^p \d x$ under the weak convergence of potentials $\{V_n\}$ in Morrey space $\mrp$, through which we prove that the functional $\lambda_1(V,m)$ is both continuous and completely continuous over $\{ (V,m) \in \mrp \times \mrp : m > 0 \text{ a.e. in } \Omega \}$. As an application, Section $4$ addresses the optimization problems for  $\lambda_1(V,m)$ over a weakly compact subset $\mathcal{W}_C \subset \{ (V,m) \in \mrp \times \mrp : m > 0 \text{ a.e. in } \Omega \}$. Finally, in Section $5$, we derive the derivative formulae of $\lambda_1(V,m)$ with respect to regular variations of $(V,m)$ in a closed subspace $M^r_0(p;\Omega)\times M^r_0(p;\Omega)$.

\section{Preliminaries}\label{sec-prem}
In this section, we outline some basic notations, definitions, and prerequisite results for the weighted eigenvalue problem \eqref{EVP}. Unless otherwise specified, we adopt the following notations, which will be used throughout the paper:
    \begin{itemize}
        \item For $y \in \mathbb{R}^N$, $B_s(y)$ denotes the open ball centered at $y$ with radius $s$.
        \item For a measurable set $E\subset \mathbb{R}^N$, $|E|$ represents the Lebesgue measure of $E$.
        \item The notation $r'$ denotes the H\"older conjugate  of $r$ i.e., $r' = \frac{r}{r-1}$.
        \item The symbol $C(\alpha,\beta,\cdots)$ used to denote the positive constant depending on $\alpha, \beta, \cdots$, which may vary line by line.
        \item The norm on the Lebesgue space $L^r(\Omega)$ ($1\leq r \leq \infty$) is denoted by $\|\cdot\|_r$.
        \item The space $\Wkp$ denotes the standard Sobolev space with zero boundary equipped with the standard Sobolev norm $\|\cdot\|_{1,p}$. 
    \end{itemize} 
We first recall the formal definition of Morrey spaces and discuss their fundamental properties. 
\begin{df}[Morrey space]\label{Morry}
  Let $r \in [1,\infty]$.  The Morrey space $M^r(\Omega)$ is defined as  $$M^r(\Omega) := \{ f \in L^1(\Omega) : \|f\|_{M^r(\Omega)} < \infty\},$$ 
    where $$
	\| f \|_{M^r(\Omega)} :=
	\sup_{\substack{y \in \Omega \\ 0 < s < \operatorname{diam}(\Omega)}}
	\frac{1}{s^{\,N/r'}} \int_{\Omega \cap B_s(y)} |f(y)| ~dy.
	$$
\end{df}
\begin{rk}\label{Remark}
\begin{itemize}
    \item [(a)] An application of H\"older inequality yields the following strict inclusions:
		$$ L^r(\Omega) \subsetneq M^r(\Omega) \subsetneq L^1(\Omega),~~ \text{for every}~~ r \in (1,\infty).$$
        \item [(b)] The Morrey space $M^r(\Omega)$ is neither reflexive nor separable (see Burenkov et al. \cite{burenkov2019}). 
        \item[(c)] The space $C_c^\infty(\Omega)$, space of all compactly supported infinitely differentiable functions in $\Omega$, is not dense in $M^r(\Omega)$ with respect to the Morrey norm. (see Burenkov et al. \cite{burenkov2019})
\end{itemize}
    \end{rk}
 \noindent For a detailed account of the regularity theory for elliptic equations involving coefficients in Morrey spaces, one can refer to the monograph of Mal\'y and Ziemer \cite{Ziemer97}, and the book of Gilbarg and Trudinger \cite{gilbarg1977}.
    
Next, we introduce a special class of Morrey spaces $M^r(p;\Omega)$ depending on the exponent $1<p<N$, which contains both the sign-changing potential $V$ and the weight  $m$ of the problem \eqref{EVP}.
\begin{df}\label{splorrey}
	Let $1<p<N$, and define the special Morrey space $M^r(p;\Omega)$ as follows:
	\begin{align*}
		M^r(p;\Omega) :=  M^r(\Omega) ; ~\text{for some}~ r > N/p.
	\end{align*}	
\end{df}
At this point, we state the following fundamental inequality for $1<p<N$ due to Morrey, subsequently generalized by Adams (see Gilbarg and Trudinger \cite{gilbarg1977}, Mal\'y and Ziemer \cite{Ziemer97}). For the general case $1<p<\infty$, the reader can refer to Pinchover and Psaradakis \cite[Theorem 2.4]{Pinchover16}.
\begin{thm}\label{M-A-THM}(Morrey-Adams theorem) Let $1<p<N$, and $V \in M^r(p;\Omega)$.  Then  there exists a constant $C(N, p, r) > 0$ such that for any $ \delta > 0 $ 
        \begin{equation}
			\int_{\Omega} |V||u|^p\d x \leq \delta \| \nabla u\|^p_{p} + \frac{C(N,p,r)}{\delta^{N/(pr-N)}} \|V\|^{pr/(pr-N)}_{M^r(p;\Omega)} \|u\|^p_{p}, \quad \text{ for all } u \in \Wkp.
		\end{equation}
\end{thm}
The Morrey-Adams Theorem \ref{M-A-THM} together with H\"older inequality ensures the validity of the weak formulation to problem \eqref{EVP} involving the potential $V$ and weight $m$.
\begin{df}\label{ws}
Let $1<p<N$,  and $A, V$ and $m$ satisfies the Assumption \ref{assump}.  Then for $\lambda\in \mathbb{R}$, a function $u\in W_0^{1,p}(\Omega)$ is a weak solution to the problem \eqref{EVP} if 
\begin{align}\label{Weak}
		\int_{\Omega} |\nabla u|^{p-2}_{A} A \nabla u \cdot \nabla v \d x + \int_{\Omega} V |u|^{p-2} u v \d x -\lambda  \int_{\Omega} m|u|^{p-2}u  v \d x =0, \,\,\forall\, v  \in W_0^{1,p}(\Omega).
	\end{align}
\end{df}
\noindent Note that any nonnegative solution to \eqref{EVP} obeys the local Harnack inequality, and consequently the solutions are locally H\"older continuous. For the proof, the reader can refer to Pinchover and Psaradakis \cite[Theorem 2.7]{Pinchover16}.

We now recall the following crucial lemma from Gupta et al. \cite[Lemma 3.2]{Gupta2026}, which will be used throughout this paper.  
\begin{lem}\label{bdd-lema}
		Let $1<p<N$. Suppose $A$, $V$,  $m$ satisfy Assumptions \ref{assump}, and $B  \subset M^r(p;\Omega) $ is a bounded subset. Then for every $V\in B$ and all $u\in \Wkp$, we have
		\begin{align}\label{bdd-1}
			\|u\|^p_{1,p} \leq C_1 \left( \int_{\Omega} |\nabla u|^{p}_A\d x + \int_{\Omega} V |u|^p\d x\right)   + C_2 \int_{\Omega} m |u|^p \d x, 
		\end{align}
        for some positive constants $C_1$ and $C_2$.
	\end{lem}
Next, we define the eigenvalues and associated eigenfunctions of \eqref{EVP} as follows:
\begin{df}[Eigenvalue]
\begin{itemize}
    \item[(i)] We say $\lambda \in \mathbb{R}$ is an eigenvalue associated with an eigenfunction $u$ if $u \in W_0^{1,p}(\Omega) \setminus \{0\}$ is a weak solution of \eqref{EVP}.
    \item[(ii)] An eigenvalue of \eqref{EVP} is said to be a principal eigenvalue if it admits a non-negative eigenfunction.
\end{itemize}  
\end{df}
\noindent Finally, we conclude this section by recalling the existence of the principal eigenvalue and its qualitative properties. 
\begin{thm}\label{ERQ}
	Let $1<p<N$, and $A$, $V$, and $m$ satisfy Assumptions \ref{assump}. Then, the Rayleigh quotient $\mathcal{R}$, defined by 
    \begin{align}
		\mathcal{R} &:= \inf_{ u \in W_0^{1,p}(\Omega)\setminus\{0\}} \frac{\int_{\Omega} |\nabla u|^{p}_A\d x + \int_{\Omega} V |u|^p\d x}{\int_{\Omega} m|u|^p\d x}\label{RQ} 
        \end{align}
 is the smallest eigenvalue for the problem \eqref{EVP}. Moreover, it is the principal eigenvalue.
\end{thm}
\begin{proof} See Gupta et al. \cite[Theorem 3.4 and 3.5]{Gupta2026}.
\end{proof}
\begin{rk}
     In view of the homogeneity of the Rayleigh quotient $\mathcal{R}$, we obtain the equivalent normalized form:
        \begin{align}
       		\mathcal{R} &= \inf_{ u \in W_0^{1,p}(\Omega)\setminus\{0\}} \left\{ \int_{\Omega} |\nabla u|^{p}_A\d x + \int_{\Omega} V |u|^p\d x : \int_{\Omega}m |u|^p\d x = 1\right\}  \label{new RQ}.
	\end{align}
\end{rk}
\noindent We denote the principal eigenvalue of \eqref{EVP} by $\lambda_1(V, m)$, i.e., $\lambda_1(V, m):= \mathcal{R}$. The following qualitative properties of the principal eigenvalue $\lambda_1(V,m)$ to \eqref{EVP} were established by Gupta et al. \cite{Gupta2026}.
\begin{thm}\label{positive}  Let $1<p<N$, $A$, $V$, and $m$ satisfy Assumptions \ref{assump}. Then the principal eigenvalue $\lambda_1(V,m)$ is unique, simple, and isolated. Moreover,  $\lambda_1(V,m)$ is  positive if $V$ satisfies the additional condition:
        \begin{equation}\label{assump-1}
            \|V^{-}\|_{M^r(p;\Omega)} \leq  \frac{\delta_0}{\left[C(N,p,r)S\right]^{(pr-N)/pr}}, ~~\text{for some}~~ 0 < \delta_0  < \frac{\Theta^p_{\Omega}}{2},
        \end{equation}
where $S$ is the Sobolev constant and $C(N,p,r)$ is the constant as in the Morrey-Adams Theorem \ref{M-A-THM}.
\end{thm}
\noindent For other properties, such as the eigenspectrum and the existence of a monotonically increasing sequence of eigenvalues to the problem \eqref{EVP}, we refer the reader to the work of Gupta et al. \cite{Gupta2026}. In the same paper, the authors proved the following boundedness result of the eigenfunction $u$ associated with the eigenvalue $\lambda(V,m)$ (see Gupta et al. \cite[Theorem 5.2]{Gupta2026}), which will be used in our present study.    
\begin{thm}\label{bounded-theorem}
	    Let $1<p<N$, and $A$, $V$, and $m$ satisfy Assumptions \ref{assump}. Suppose  $u$ is an eigenfunction associated with the eigenvalue $\lambda(V,m)$ of the problem \eqref{EVP}. Then $u \in L^\infty(\Omega)$.
	\end{thm}
    
\section{Properties of the functional $\lambda_1(V,m)$}
This section is devoted to studying various properties of the functional $\lambda_1(V,m)$ for $(V,m) \in \mrp \times \mrp$ with $m > 0$ a.e. in $\Omega$. We start with the following concavity and convexity of  $\lambda_1(V,m)$ with respect to potential $V$ and weight $m$, respectively. 
\begin{prop}
    The functional $\lambda_1(V,m)$ is concave with respect to the potential $V$.
\end{prop}
\begin{proof}
Let $ t \in [0,1]$, and $tV_1 + (1-t)V_2 \in M^r(p;\Omega)$ for $V_1, V_2 \in \mrp$. By the definition of $\lambda_1(V,m)$, we obtain 
    \begin{align}\label{conc-1}
        t\lambda_1(V_1,m) +& (1-t) \lambda_1(V_2,m)\nonumber \\ & \leq t \left( \int_{\Omega} |\nabla u|^p_A \d x + \int_{\Omega} V_1|u|^p \d x \right)+ (1-t)\left( \int_{\Omega} |\nabla u|^p_A \d x + \int_{\Omega} V_2|u|^p \d x \right)\nonumber\\
        & = \int_{\Omega} |\nabla u|^p_A \d x + \int_{\Omega} (tV_1 + (1-t)V_2) |u|^p \d x,
    \end{align}
for all $u \in \Wkp$ such that $\int_{\Omega}m |u|^p\d x = 1$. Now taking the infimum over set $\{ u \in \Wkp : \int_{\Omega}m |u|^p\d x = 1\}$, we conclude that 
\begin{align*}
    t\lambda_1(V_1,m) + (1-t)\lambda_1(V_2,m)  \leq \lambda_1 (tV_1 + (1-t)V_2,m).
\end{align*}
This shows that $\lambda_1(V,m)$ is concave with respect to the potential $V$.
\end{proof}
\begin{prop} Suppose that the potential $V$ satisfies the additional  assumption \eqref{assump-1}. Then
\begin{itemize}
    \item[(i)] the functional $\frac{1}{\lambda_1(V,m)}$ is convex with respect to the weight $m$.
    \item[(ii)] the functional $\frac{1}{\lambda_1(V,m)}$ is strictly convex with respect to $m$ on the set $G_K := \{ m \in \mrp : m > 0 \text{ a.e. in } \Omega \text{ and } \int_{\Omega} m \,\d x = K\}$.
\end{itemize}
\end{prop}
\begin{proof}
(i) By the Theorem \ref{positive}, we have $\lambda_1(V,m) > 0$. Then, we have  $$\frac{1}{\lambda_1(V,m)}= \sup_{ u \in \Wkp \setminus \{0\}} \frac{\int_{\Omega} m|u|^p \d x}{\int_{\Omega} |\nabla u|^p_A \d x + \int_{\Omega} V|u|^p \d x}.$$
Let $ t \in [0,1]$, and $m_1, m_2 \in \mrp$ with both $m_1,m_2 >0$ a.e. in $\Omega$. Then $ m^* := tm_1 + (1-t)m_2 \in M^r(p;\Omega)$ with $m^*>0$ a.e. in $\Omega$. Thus
\begin{align*}
\frac{1}{\lambda_1(V,m^*)}&= \sup_{ u \in \Wkp \setminus \{0\}}\frac{\int_{\Omega} m^*|u|^p \d x }{\int_{\Omega} |\nabla u|^p_A \d x + \int_{\Omega} V|u|^p \d x}\\
 &= \sup_{ u \in \Wkp \setminus \{0\}} \frac{\int_{\Omega} t m_1|u|^p \d x + \int_{\Omega} (1-t) m_2|u|^p \d x}{\int_{\Omega} |\nabla u|^p_A \d x + \int_{\Omega} V|u|^p \d x}\\
     & =\sup_{ u \in \Wkp \setminus \{0\}} \left[ \frac{t\int_{\Omega}  m_1|u|^p \d x}{\int_{\Omega} |\nabla u|^p_A \d x + \int_{\Omega} V|u|^2 \d x} +  \frac{(1-t)\int_{\Omega}  m_2|u|^p \d x}{\int_{\Omega} |\nabla u|^p_A \d x + \int_{\Omega} V|u|^p \d x}\right]\\
     % & \leq \sup_{0 \neq u \in \Wkp}\frac{\int_{\Omega} t m_1|u|^p \d x}{\int_{\Omega} |\nabla u|^p_A \d x + \int_{\Omega} V|u|^2 \d x}\\
     % & \hspace{5cm}+\sup_{0 \neq u \in \Wkp}\frac{\int_{\Omega} (1-t) m_2|u|^p \d x}{\int_{\Omega} |\nabla u|^p_A \d x + \int_{\Omega} V|u|^2 \d x}\\
 & \leq \frac{t}{\lambda_1(V,m_1)} + \frac{(1-t)}{\lambda_1(V,m_2)}. 
\end{align*}
This implies that the functional $\frac{1}{\lambda_1(V,m)}$ is convex with respect to $m$.\\
(ii) Suppose that $m_1, m_2 \in G_K$. Let $u^*, u_1$ and $u_2$ be positive eigenfunctions associated with $\lambda_1(V,m^*)$, $\lambda_1(V,m_1)$ and $\lambda_1(V,m_2)$, respectively. Now for $0 \leq  t \leq 1 $, and if the equality holds in the above, then we get
\begin{align*}
    \frac{\int_{\Omega} m_1|u^*|^p \d x}{\int_{\Omega} |\nabla u^*|^p_A \d x + \int_{\Omega} V|u^*|^p \d x} = \frac{1}{\lambda_1(V,m_1)},
%\end{align*}
\text{ and } 
%\begin{align*}
     \frac{\int_{\Omega} m_2|u^*|^p \d x}{\int_{\Omega} |\nabla u^*|^p_A \d x + \int_{\Omega} V|u^*|^p \d x} = \frac{1}{\lambda_1(V,m_2)}.
\end{align*}
As a consequence, we have 
$$ -\Delta_{p,A}(u^*) + V (u^*)^{p-1} = \lambda_1(V,m_1)m_1(u^*)^{p-1},$$
and 
$$ -\Delta_{p,A}(u^*) + V (u^*)^{p-1} = \lambda_1(V,m_2)m_2(u^*)^{p-1}.$$
Thus, from the weak formulation, we have
$$ \int_{\Omega}\left[ \lambda_1(V,m_1) m_1 - \lambda_1(V,m_2) m_2\right] (u^*)^{p-1}\varphi \d x = 0, \text{ for all } \varphi \in \Wkp.$$
Since $u^* > 0$, hence we get $\lambda_1(V,m_1) m_1 = \lambda_1(V,m_2) m_2$. Next, integrating over $\Omega$ and recalling that $m_1,m_2 \in G_K$, we find  $\lambda_1(V,m_1) = \lambda_1(V,m_2)$. This completes the proof.
\end{proof}
Next, we discuss the continuity property of the functional $\lambda_1(V,m)$.
\begin{thm}\label{continuous-thm} 
The functional $\lambda_1(V,m)$ is continuous on  $\{(V,m) \in \mrp \times \mrp : m>0 \text{ a.e. in } \Omega\}.$ 
\end{thm} 
\begin{proof}
    Let $ (V_n,m_n) \to (V,m)$ in $\mrp  \times \mrp$ as $n \to \infty$, where $m_n > 0$ a.e. in $\Omega$.  Let $u_n \in \Wkp$ be a positive  eigenfunction associated with  $\lambda_1(V_n,m_n)$ for every $n \in \mathbb{N}$. Without loss of generality, let us assume that $\int_{\Omega} m_n|u_n|^p \d x = 1$  for each $n \in \mathbb{N}$. Then, for every $n$ we have 
    \begin{align}\label{conti-1}
    \int_{\Omega} |\nabla u_n|^p_A \d x + \int_{\Omega} V_n|u_n|^p \d x = \lambda_1(V_n,m_n). 
    \end{align}
    Let $u$ be a positive eigenfunction corresponding to  $\lambda_1(V,m)$ such that $\int_{\Omega}m|u|^p \d x =1$.
   Since $V_n \to V$ in $\mrp$, hence (up to a subsequence) $V_n \to V$ a.e. in $\Omega$. Thus, $V_n|u|^p \to V|u|^p$ a.e. in $\Omega$, and $|V_n|u|^p| \leq C |V|$ (as $u \in L^{\infty}(\Omega)$ by  Theorem \ref{bounded-theorem}), where $C$ is a positive constant independent of $n$. Since $|V| \in \mrp \subset L^1(\Omega)$, it follows from the dominated convergence theorem that
   \begin{align}\label{conti-2}
      \int_{\Omega}V_n|u|^p \d x \to \int_{\Omega}V|u|^p \d x,\quad \text{as } n \to \infty.
   \end{align}
Similarly, as $m_n\rightarrow m $ in $\mrp$, we conclude that 
\begin{align}\label{conti-3}
      \int_{\Omega}m_n|u|^p \d x \to \int_{\Omega}m|u|^p \d x,\quad \text{as } n \to \infty.
   \end{align}
This shows that $ \int_{\Omega} m_n |u|^p \d x  \to 1,$ as  $ n \to \infty$. 
Choosing $v_n := \frac{u}{(\int_{\Omega}m_n|u|^p \d x)^{\frac{1}{p}}}$ as an admissible function in \eqref{RQ} for $\lambda_1(V_n,m_n)$, and then applying Morrey-Adams Theorem \ref{M-A-THM},  we have 
\begin{align}\label{conntti-1}
    \lambda_1(V_n,m_n)& \leq \int_{\Omega}|\nabla v_n|^p_A \d x + \int_{\Omega} V_n |v_n|^p \d x \nonumber \\
    &= \frac{\int_{\Omega} |\nabla u|^p_A \d x + \int_{\Omega} V_n|u|^p \d x}{\int_{\Omega}m_n|u|^p \d x}  \\
   % &\leq \frac{\int_{\Omega} |\nabla u|^p_A \d x + \int_{\Omega} |V_n||u|^p \d x}{\int_{\Omega}m_n|u|^p \d x} \nonumber \\
    &\leq \frac{   \Theta_{\Omega}^{-p} \|u\|_{1,p}^p+\delta \| \nabla u\|^p_{p} + \frac{C(N,p,r)}{\delta^{N/(pr-N)}} \|V_n\|^{pr/(pr-N)}_{M^r(p;\Omega)} \|u\|^p_{p}}{\int_{\Omega}m_n|u|^p \d x}.\nonumber
\end{align}
By using the boundedness of sequences $\{V_n\}$  in $ \mrp$ and $\left\{\frac{1}{\int_{\Omega}m_n|u|^p \d x} \right\}$ in $\mathbb{R}$, we get  $\lambda_1(V_n,m_n) < C $, for some constant $C > 0$ independent of $n$. Then, applying the Lemma \ref{bdd-lema}, we have 
$$\|u_n\|^p_{1,p} \leq C_1  \lambda_1(V_n,m_n) + C_2 \int_{\Omega}m_n|u_n|^p \d x \leq C_1 C + C_2.$$
This implies that $\{u_n\} $ is bounded in $\Wkp$, and hence there exists  $w \in \Wkp$ such that $u_n \to w$ in $\Lp$ and $u_n \to w$ a.e. in $\Omega$ (up to a subsequence). The boundedness of sequences $\{V_n\}$ in $ \mrp$,  $\{u_n\}$ in $ \Wkp$, and together with the Morrey-Adamas Theorem \ref{M-A-THM}, implies that 
    \begin{align}\label{conti-5}
        \int_{\Omega} |V_n||u_n|^p \d x \leq C_3,
    \end{align}
    where $C_3 > 0 $ is independent of $n$. Moreover,
    \begin{align}\label{conti-6}
        \int_{\Omega} |V_n||w|^p \d x \leq C_4,
    \end{align} 
    for some constant $C_4 > 0 $ independent of $n$.

Next, we claim $\int_{\Omega}V_n|u_n|^p \d x \to \int_{\Omega}V|w|^p \d x$ as $ n \to \infty$. Notice that 
\begin{align}\label{conti-7}
    \left| \int_{\Omega} V_n|u_n|^p \d x - \int_{\Omega}V|w|^p \d x   \right| \leq \int_{\Omega} |V_n|\left| |u_n|^p-|w|^p\right| \d x + \int_{\Omega} |V_n-V||w|^p \d x.
\end{align}
 By virtue of the inequality $|c^p - \d^p| \leq p|c-d|(c^{p-1} + d^{p-1}), ~\text{for}~~c, d \geq 0,$ and then invoking the H\"older inequality, we  obtain
        \begin{align}\label{conti-8}
        &\int_{\Omega}|V_n| \left|  |u_n|^p - |w|^p\right|\d x \leq p \int_{\Omega} |V_n| |u_n - w| \left| |u_n|^{p-1} + |w|^{p-1} \right|\d x \nonumber\\
       % & \hspace{1cm}= p \int_{\Omega} |V|^{1/p} |u_n - u| |V|^{1/p'} \left| |u_n|^{p-1} + |w|^{p-1} \right| \d x\nonumber \\
        & \hspace{1cm}\leq p \left( \int_{\Omega} |V_n| |u_n - w|^p \d x\right)^{1/p} \left(\int_{\Omega}  |V_n|  \left| |u_n|^{p-1}+ |w|^{p-1} \right|^{\frac{p}{p-1}} \d x\right)^{1/p'} \nonumber\\
        &\hspace{1cm}\leq C(p)\left( \int_{\Omega} |V_n| |u_n - w|^p\d x\right)^{1/p} \left(\int_{\Omega}  |V_n|   |u_n|^p\d x + \int_{\Omega}  |V_n|  |w|^p \d x \right)^{1/p'}.
        \end{align}
       Let $\varepsilon > 0$  be arbitrary. By the $L^p$-convergence of the sequence $\{u_n\}$, we can find $N_{\varepsilon} \in \mathbb{N}$ satisfying
       \begin{align}\label{conti-9}
           \int_{\Omega} |u_n-w|^p \d x \leq \varepsilon^{\frac{pr}{pr-N}}, \text{ for all } n \geq N_{\varepsilon}.
       \end{align} 
Then applying the Morrey-Adams Theorem \ref{M-A-THM} for $\varepsilon >0 $, and  using the boundedness of sequence $\{u_n\}$ in $\Wkp$ together with inequality \eqref{conti-9}, we deduce that
	\begin{align}\label{conti-10}
			\left( \int_{\Omega} \!\!|V_n| |u_n - w|^p\d x\right)^{\frac{1}{p}} \!\!&\leq  \left(\varepsilon \!\int_\Omega |\nabla (u_n - w)|^p  \d x +   \frac{ C(N,p,r)}{\varepsilon^{N/(pr-N)}}\|V_n\|^{\frac{pr}{pr-N}}_{M^r(p;\Omega)}\!\!\int_\Omega |u_n - w|^p \d x  \right)^{\frac{1}{p}} \nonumber\\
            & \leq C_5 \varepsilon, 
		\end{align}
        for some constant $C_5 > 0$, and all $n > N_{\varepsilon}$. 
Combining the estimates \eqref{conti-5},\eqref{conti-6} and \eqref{conti-10} in \eqref{conti-8}, we get for all $n > N_{\varepsilon}$
        \begin{align}\label{conti-11}
            \int_{\Omega}|V_n| \left|  |u_n|^p - |w|^p\right|\d x \leq C_6 \varepsilon,
        \end{align}
        for some  constant $C_6>0$. Next, arguing similar as in \eqref{conti-2}, we also have
 \begin{align}\label{conti-12}
     \int_{\Omega} |V_n-V||w|^p \d x \to 0 \text{ as } n \to \infty.
 \end{align}
Therefore, using \eqref{conti-11} and \eqref{conti-12} in \eqref{conti-7}, we obtain
% \begin{align}
%     \left| \int_{\Omega} V_n|u_n|^p \d x - \int_{\Omega}V|w|^p \d x   \right| \leq C_7 \varepsilon \text{ for all } n > N_{\varepsilon},
% \end{align}
% where $C_7$ is positive constant independent of $n$. Hence, we have
    \begin{align}\label{conti-13}
        \int_{\Omega} V_n |u_n|^p \d x \to \int_{\Omega} V|w|^p \d x \text{ as }  n \to \infty.
    \end{align}
Similarly, proceeding as above,  we have
\begin{align}\label{conti-14}
    \int_{\Omega} m_n|u_n|^p \d x \to \int_{\Omega}m|w|^p \d x \text{ as }n \to \infty.
\end{align}
This implies that $\int_{\Omega} m|w|^p \d x = 1$. Thus, by 
taking $w$ as an admissible function in \eqref{RQ} for $\lambda_1(V,m)$ and then using weak lower semicontinuity of $\int_{\Omega} |\nabla u|^p_A \d x$ together with \eqref{conti-13}, we get
    \begin{align}\label{conti-15}
        \lambda_1(V,m) \leq  \int_{\Omega} |\nabla w|^p_A \d x + \int_{\Omega} V|w|^p \d x &\leq \liminf_{n \to \infty} \left( 
     \int_{\Omega} |\nabla u_n|^p_A \d x + \int_{\Omega} V_n|u_n|^p \d x\right) \nonumber\\
     &=\liminf_{n \to \infty} \lambda_1(V_n,m_n). 
    \end{align}
    Furthermore, from \eqref{conntti-1}, we deduce that
    \begin{align}\label{conntti-2}
 \limsup_{n \to \infty} \lambda_1(V_n,m_n) \leq \lambda_1(V,m).
    \end{align}
Hence, from \eqref{conti-15} and \eqref{conntti-2}  the result follows.
    % Taking the limit as $ n \to \infty$ in \eqref{conti-1} and using \eqref{conti-3}, \eqref{conti-4} and \eqref{conti-15}, we obtain the following:
    % \begin{align*}
    %     \lambda_1(V,m) \leq \liminf_{n \to \infty} 
    %  \lambda_1(V_n,m_n) \leq  \limsup_{n \to \infty} \lambda_1(V_n,m_n) \leq \lambda_1(V,m).
    % \end{align*}
    % Hence, $ \lambda_1(V_n,m_n) \to \lambda_1(V,m)$ for $V_n \to V$ and $m_n \to m$ as $n \to \infty$. Thus $\lambda_1(V,m)$ is continuous.
\end{proof}
At this stage,  we further show that the functional $\lambda_1(V,m)$ is completely continuous, in the sense that if $V_n \rightharpoonup V$ and $ m_n \rightharpoonup m$ in $\mrp$ (with $m_n, m >0$ a.e. in $\Omega$) implies  $\lambda_1(V_n,m_n) \to \lambda_1(V,m)$ as $n \to \infty$. To establish this result, we now prove several lemmas concerning the convergence of integrals of the form $\int_\Omega V_n |u_n|^p \d x$ under the weak convergence of the sequence of potentials $\{V_n\}$ in $\mrp$. Furthermore, these results play a crucial role in analyzing optimization problem and the differentiability property of the functional $\lambda_1(V, m)$. 
\begin{lem}\label{convergence-1}
    Let  $u$ be a positive eigenfunction associated with the eigenvalue  $\lambda_1(V,m)$. If $V_n \rightharpoonup V$  in $\mrp$, then 
    $$\int_{\Omega}V_n|u|^p \d x \to \int_{\Omega}V|u|^p \d x,\quad \text{as } n \to \infty.$$
\end{lem}
\begin{proof}
   Since $u$ is an eigenfunction corresponding to $\lambda_1(V,m)$, Theorem \ref{bounded-theorem} implies that $u \in L^\infty(\Omega)$. Fixing $u$, define the map $F : \mrp \to \mathbb{R}$ by $\displaystyle{F(V) = \int_{\Omega}V|u|^p \d x}$.  Then, using the continuous embedding $\mrp \subset L^1(\Omega)$, it follows that $F$ is a bounded linear functional. Hence, by the notion of weak convergence, the claim follows.
   \end{proof}
\begin{lem}\label{convergence-2}
   If $V_n \rightharpoonup V$  in $\mrp$, then for every $u \in \Wkp$, we have
   $$\int_{\Omega} V_n|u|^p \d x \to \int_{\Omega}V|u|^p \d x \text{ as } n \to \infty.$$
\end{lem}
\begin{proof}
  Let $\varepsilon > 0$ be arbitrary. Since $u \in \Wkp$, and $C_c^\infty(\Omega) $ is dense in $\Wkp$, hence there exists $\varphi_\varepsilon \in C_c^\infty(\Omega)$  such that
    \begin{align}\label{wk-00}
         \|u-\varphi_\varepsilon\|_{1,p} < \varepsilon^{\frac{r}{pr-N}}.
    \end{align}
    Now
    \begin{align}\label{wk-01}
        & \left| \int_{\Omega}\!\! V_n|u|^p \d x - \int_{\Omega}\!\!V|u|^p \d x   \right| =           \left| \int_{\Omega} \!\!\left(V_n|u|^p - V_n|\varphi_\varepsilon|^p + V_n|\varphi_\varepsilon|^p- V|\varphi_\varepsilon|^p + V|\varphi_\varepsilon|^p  - V|u|^p  \right) \d x   \right| \nonumber\\
         &\hspace{2cm} \leq \int_{\Omega}|V_n|||u|^p-|\varphi_\varepsilon|^p| \d x + \int_{\Omega}|V_n-V||\varphi_\varepsilon|^p \d x + \int_{\Omega} |V|||\varphi_\varepsilon|^p-|u|^p|\d x.
    \end{align}
  By using the inequality $|c^p - \d^p| \leq p|c-d|(c^{p-1} + d^{p-1}), ~\text{for}~~c, d \geq 0,$ and then invoking the H\"older inequality, we obtain
        \begin{align}\label{wk-02}
        &\int_{\Omega}|V_n| \left|  |u|^p - |\varphi_\varepsilon|^p\right|\d x \leq p \int_{\Omega} |V_n| |u - \varphi_\varepsilon| \left| |u|^{p-1} + |\varphi_\varepsilon|^{p-1} \right|\d x \nonumber\\
        & \hspace{1cm}\leq p \left( \int_{\Omega} |V_n| |u - \varphi_\varepsilon|^p \d x\right)^{1/p} \left(\int_{\Omega}  |V_n|  \left| |u|^{p-1}+ |\varphi_\varepsilon|^{p-1} \right|^{\frac{p}{p-1}} \d x\right)^{1/p'} \nonumber\\
        &\hspace{1cm}\leq C(p)\left( \int_{\Omega} |V_n| |u - \varphi_\varepsilon|^p\d x\right)^{1/p} \left(\int_{\Omega}  |V_n|   |u|^p\d x + \int_{\Omega}  |V_n|  |\varphi_\varepsilon|^p \d x \right)^{1/p'}.
        \end{align}
        Since $V_n \rightharpoonup V$ in $\mrp$ implies that
        $\{V_n\}$ is bounded in  $ \mrp$. Thus, invoking  the Morrey-Adamas Theorem \ref{M-A-THM}, implies that 
    \begin{align}\label{wk-03}
        \int_{\Omega} |V_n||u|^p \d x \leq C_1,
    \end{align}
    where $C_1 > 0 $ is  constant independent of $n$. As $\varphi_\varepsilon \in C_c^\infty(\Omega)$, we also have
    \begin{align}\label{wk-04}
        \int_{\Omega} |V_n||\varphi_\varepsilon|^p \d x \leq C_2,
    \end{align} 
    for some constant $C_2 > 0 $ independent of $n$. Further, by Poincar\'e inequality, we obtain
    \begin{align}\label{wk-004}
        \|u-\varphi_\varepsilon\|_p^p \leq C_3 \|u-\varphi_\varepsilon\|_{1,p}^p \leq C_3 \varepsilon^{\frac{pr}{pr-N}}.
    \end{align} 
    Again using the Morrey-Adamas Theorem \ref{M-A-THM} for $\delta = \varepsilon$ and \eqref{wk-004}, we have
    \begin{align}\label{wk-05}
			\int_{\Omega}|V_n| |u - \varphi_\varepsilon|^p\d x &\leq  \varepsilon \|u-\varphi_\varepsilon\|_{1,p}^p +   \frac{ C(N,p,r)}{\varepsilon^{N/(pr-N)}}\|V_n\|^{\frac{pr}{pr-N}}_{M^r(p;\Omega)}   \|u-\varphi_\varepsilon\|_{p}^p \nonumber\\
            & \leq C_4 \varepsilon, 
		\end{align}
        for some constant $C_4> 0$ independent of $n$. Here, in the above, we have used the boundedness of $\{V_n\}$ in $\mrp$.
        Then, combining \eqref{wk-03},\eqref{wk-04} and \eqref{wk-05} in  \eqref{wk-02}, we get
        \begin{align}\label{wk-06}
            \int_{\Omega}|V_n| \left|  |u|^p - |\varphi_\varepsilon|^p\right|\d x  \leq C_5 \varepsilon,
        \end{align}
        for some constant $C_5 >0$ independent of $n$. Also, repeating  the same argument as above, we obtain 
        \begin{align}\label{wk-07}
            \int_{\Omega}|V| \left|   |\varphi_\varepsilon|^p  - |u|^p  \right|\d x \leq C_6 \varepsilon,
        \end{align}
        for some $C_6 >0$ constant. Since $V_n \rightharpoonup V$ in $\mrp$ and $\varphi_\varepsilon \in C_c^\infty(\Omega)$, hence the notion of weak convergence implies
        \begin{align}\label{wk-08}
            \int_{\Omega}|V_n-V||\varphi_\varepsilon|^p \d x \to 0 \text{ as } n \to \infty.
        \end{align}
        Therefore, combining \eqref{wk-06},\eqref{wk-07} and \eqref{wk-08} in \eqref{wk-01}, we deduce that
        $$ \int_{\Omega} V_n|u|^p \d x \to \int_{\Omega}V|u|^p \d x, $$ as $n \to \infty.$
\end{proof}

\begin{lem}\label{convergence-3}
 If $V_n \rightharpoonup V$ in $\mrp$, and  $u_n \rightharpoonup  w $ in $\Wkp$, then 
 $$\int_{\Omega}V_n|u_n|^p \d x \to \int_{\Omega}V|w|^p \d x,\quad \text{as } n \to \infty.$$
\end{lem}
\begin{proof}
 Since $u_n \rightharpoonup w$  in $ \Wkp$, hence $\{u_n\}$ bounded in  $ \Wkp$. Then,  following the similar arguments as in the Lemma \ref{convergence-2}, we have
    \begin{align}\label{wk-cont-5}
        \int_{\Omega} |V_n||u_n|^p \d x \leq C_1 \text{ and } \int_{\Omega} |V_n||w|^p \d x \leq C_2,
    \end{align}
    for some constants $C_1, C_2 >0$, which are independent of $n$. Observe that 
\begin{align}\label{wk-cont-7}
    \left| \int_{\Omega} V_n|u_n|^p \d x - \int_{\Omega}V|w|^p \d x   \right| \leq    \int_{\Omega} |V_n|\left| |u_n|^p-|w|^p\right| \d x + \left|\int_{\Omega} (V_n-V)|w|^p \d x \right|.
\end{align}
 By virtue of the inequality $|c^p - \d^p| \leq p|c-d|(c^{p-1} + d^{p-1}), ~\text{for}~~c, d \geq 0,$ and then using the H\"older inequality, we  obtain
        \begin{align}\label{wk-cont-08}
        &\int_{\Omega}|V_n| \left|  |u_n|^p - |w|^p\right|\d x \leq p \int_{\Omega} |V_n| |u_n - w| \left| |u_n|^{p-1} + |w|^{p-1} \right|\d x \nonumber\\
       % & \hspace{1cm}= p \int_{\Omega} |V|^{1/p} |u_n - u| |V|^{1/p'} \left| |u_n|^{p-1} + |w|^{p-1} \right| \d x\nonumber \\
        & \hspace{1cm}\leq p \left( \int_{\Omega} |V_n| |u_n - w|^p \d x\right)^{1/p} \left(\int_{\Omega}  |V_n|  \left| |u_n|^{p-1}+ |w|^{p-1} \right|^{\frac{p}{p-1}} \d x\right)^{1/p'} \nonumber\\
        &\hspace{1cm}\leq C(p)\left( \int_{\Omega} |V_n| |u_n - w|^p\d x\right)^{1/p} \left(\int_{\Omega}  |V_n|   |u_n|^p\d x + \int_{\Omega}  |V_n|  |w|^p \d x \right)^{1/p'}.
        \end{align}
       Further, by the compact embedding $ \Wkp \subset \Lp$, we have $u_n \to w$ in $\Lp$ as $n \to \infty$.  Let $\varepsilon > 0$  be arbitrary. Then, we can find $N_{\varepsilon} \in \mathbb{N}$ satisfying
       \begin{align}\label{wk-cont-09}
           \|u_n-w\|_{p}  \leq \varepsilon^{\frac{r}{pr-N}}, \text{ for all } n \geq N_{\varepsilon}.
       \end{align} 
Then, applying the Morrey-Adams Theorem \ref{M-A-THM} for $ \delta = \varepsilon >0 $,   using the boundedness of the sequence $\{u_n\}$ in $\Wkp$ and $\{V_n\}$ in $\mrp$ together with inequality \eqref{wk-cont-09}, we deduce that
	\begin{align}\label{wk-cont-10}
			 \int_{\Omega} \!\!|V_n| |u_n - w|^p\d x &\leq  \varepsilon \|u_n-w\|^p_{1,p} +   \frac{ C(N,p,r)}{\varepsilon^{N/(pr-N)}}\|V_n\|^{\frac{pr}{pr-N}}_{M^r(p;\Omega)} \|u_n-w\|_p^p  \nonumber\\
            & \leq C_3 \varepsilon, 
		\end{align}
        for some constant $C_3 > 0$ independent of $n$ and all $n \geq N_{\varepsilon}$. 
Combining estimates \eqref{wk-cont-5} and \eqref{wk-cont-10} in \eqref{wk-cont-08}, we get  for all $n \geq N_{\varepsilon}$
        \begin{align}\label{wk-cont-11}
            \int_{\Omega}|V| \left|  |u_n|^p - |w|^p\right|\d x \leq C_4\varepsilon,
        \end{align}
        for some  constant $C_4>0$ independent of $n$.  Next, using  Lemma \ref{convergence-2}, we have
 \begin{align}\label{wk-cont-12}
     \left| \int_{\Omega} (V_n-V)|w|^p \d x \right| \to 0 \text{ as } n \to \infty.
 \end{align}
Therefore, using \eqref{wk-cont-11} and \eqref{wk-cont-12} in \eqref{wk-cont-7}, we obtain
    \begin{align*}
        \int_{\Omega} V_n |u_n|^p \d x \to \int_{\Omega} V|w|^p \d x,  
    \end{align*}
    as $ n \to \infty$.
\end{proof}
Consequently,  we arrive at the following result:
\begin{thm}\label{Lemma-weak-continuous} The functional $\lambda_1(V,m)$ is completely continuous on $\{ (V, m)\in \mrp \times \mrp : m > 0 \text{ a.e.  in } \Omega \}$, i.e., if $V_n \rightharpoonup V$ in $\mrp$, and $ m_n \rightharpoonup m$ in $\mrp$ with $m_n, m>0$ a.e. in $\Omega$ for all $n\in\mathbb{N}$, then $\lambda_1(V_n,m_n) \to \lambda_1(V,m)$.
\end{thm}
\begin{proof}
 Suppose that $u_n$ and $w$ are eigenfunctions corresponding to  eigenvalues $\lambda_1(V_n,m_n)$ and $\lambda_1(V,m)$, respectively such that $\displaystyle{\int_{\Omega}m_n|u_n|^p \d x = 1}$ and $\displaystyle{\int_{\Omega}m|w|^p \d x = 1}$. Then, the Lemma \ref{convergence-1} implies that  $ \displaystyle{\int_{\Omega}m_n|w|^p \d x \to \int_{\Omega}m|w|^p \d x=1 \text{ as } n \to \infty }$. Now, by the definition of $\lambda_1(V_n,m_n)$, and then applying the Morrey-Adams Theorem \ref{M-A-THM}, we get
\begin{align*}
    \lambda_1(V_n,m_n) &\leq \frac{\int_{\Omega} |\nabla w|^p_A \d x + \int_{\Omega} V_n|w|^p \d x}{\int_{\Omega}m_n|w|^p \d x} \nonumber \\
    &\leq \frac{   \Theta_{\Omega}^{-p} \|w\|_{1,p}^p+\delta \| \nabla w\|^p_{p} + \frac{C(N,p,r)}{\delta^{N/(pr-N)}} \|V_n\|^{pr/(pr-N)}_{M^r(p;\Omega)} \|w\|^p_{p}}{\int_{\Omega}m_n|w|^p \d x}.
\end{align*}
By using the boundedness of sequences $\{V_n\}$  in $ \mrp$ and $\left\{\frac{1}{\int_{\Omega}m_n|w|^p \d x} \right\}$ in $\mathbb{R}$, we get  $\lambda_1(V_n,m_n) < C $, for some constant $C > 0$ independent of $n$. Thus, from Lemma \ref{bdd-lema}, we have
$$\|u_n\|^p_{1,p} \leq C_1  \lambda_1(V_n,m_n) + C_2 \int_{\Omega}m_n|u_n|^p \d x \leq C_1 C + C_2.$$
This shows that $\{u_n\}$ is bounded in $\Wkp$. Hence there exists $u \in \Wkp$ such that $u_n \rightharpoonup u$ in $\Wkp$  and $u_n \to u$ a.e. in $\Omega$ (up to a subsequence). Hence, using the Lemmas  \ref{convergence-2} and \ref{convergence-3}, we have 
\begin{equation}\label{weak-cont-1}
\left\{
\begin{aligned}
\int_{\Omega}(V-V_n)|u_n|^p \d x \to 0  \text{ as } n \to \infty,\\
\int_{\Omega}(V_n-V)|w|^p \d x \to 0 \text{ as } n \to \infty.
\end{aligned}
\right.
\end{equation}
Next, using $u_n $ as an admissible function in the definition of $\lambda_1(V,m)$, we obtain
\begin{align}\label{weak-cont-2}
\lambda_1(V,m) \int_{\Omega}m|u_n|^p \d x &\leq \int_{\Omega} |\nabla u_n|^p_A \d x + \int_{\Omega} V|u_n|^p \d x \nonumber\\
& = \int_{\Omega} |\nabla u_n|^p_A \d x + \int_{\Omega} V|u_n|^p \d x + \int_{\Omega} V_n|u_n|^p \d x -\int_{\Omega} V_n|u_n|^p \d x \nonumber\\
& = \lambda_1(V_n,m_n) + \int_{\Omega}(V- V_n)|u_n|^p \d x.
\end{align}
Similarly, using $w$ as an admissible function in  the definition of   $\lambda_1(V_n,m_n)$, we arrive at
\begin{align}\label{weak-cont-3}
\lambda_1(V_n,m_n) \int_{\Omega}m_n|w|^p \d x \leq  \lambda_1(V,m) + \int_{\Omega}(V_n- V)|w|^p \d x.
% \implies & \lambda_1(V_n,m_n)  \leq  \frac{\lambda_1(V,m) + \int_{\Omega}(V_n- V)|w|^p \d x}{\int_{\Omega}m_n|w|^p \d x }.
\end{align}
Further, as $m_n,m \in \mrp$, we can apply Lemma  \ref{convergence-3} to obtain
\begin{equation*}\label{weak-cont-5}
\int_{\Omega}(m-m_n)|u_n|^p \d x \to 0  \text{ as } n \to \infty.
\end{equation*}
This implies that
\begin{align}\label{weak-cont-6}
     \int_{\Omega} m |u_n|^p \d x  \to 1 \text{ as } n \to \infty. 
\end{align}
By combining \eqref{weak-cont-2} and \eqref{weak-cont-3}, we have 
\begin{align}\label{weak-cont-4}
 \lambda_1(V,m) \int_{\Omega}m|u_n|^p \d x   - \int_{\Omega}(V- V_n)|u_n|^p \d x  &\leq  \lambda_1(V_n,m_n)\nonumber\\ 
 &\leq  \frac{\lambda_1(V,m) + \int_{\Omega}(V_n- V)|w|^p \d x}{\int_{\Omega}m_n|w|^p \d x }.
\end{align}
 Thus, using \eqref{weak-cont-1}  and \eqref{weak-cont-6} in \eqref{weak-cont-4}, we  get  $\lambda_1(V_n,m_n) \to \lambda_1(V,m)$ as $n \to \infty$.
\end{proof}

\section{Optimization of $\lambda_1(V,m)$ over a weakly compact subset}
Since the Morrey space $\mrp$ fails to be separable and reflexive, we study the existence of optimal potentials and positive weights for the principal eigenvalue $\lambda_1(V,m)$ over a weakly compact subset.

 Let $\mathcal{W}_C\subset \{ (V, m)\in \mrp \times \mrp : m > 0 \text{ a.e.  in } \Omega \} \subset \mrp \times \mrp$ be a weakly compact subset. Since $\mathcal{W}_C$ is bounded in $\mrp \times \mrp$ (see \cite[Proposition 2.8.1]{Robert98}), hence there exists a constant $K>0$ such that $$\|V\|_{\mrp} \leq K \text{ and } \|m\|_{\mrp} \leq K,$$ for every $(V, m)\in \mathcal{W}_C$. We now show that the set of principal eigenvalues $\lambda_1(V,m)$ is bounded from above over the weakly compact set $\mathcal{W}_C$.
\begin{lem}\label{bdd_lemma}
    There exists $C(N,p,r,S,K,\Theta_\Omega)> 0$ such that $$\lambda_1(V, m) \leq C, $$
for all $(V, m)\in \mathcal{W}_C$.
\end{lem}
\begin{proof}
Let $(V, m) \in \mathcal{W}_C$. Let $\varphi \in C_c^\infty(\Omega)$. By applying the Morrey–Adams Theorem \ref{M-A-THM} for a fix $\delta > 0$, followed by the Poincar\'e inequality, we have
\begin{align}\label{max-7}
 \int_{\Omega} |V||\varphi|^p \d x  &\leq    \delta \|\varphi\|^p_{1,p} + \frac{ C(N,p,r)}{\delta^{N/(pr-N)}}\|V \|^{\frac{pr}{pr-N}}_{M^r(p;\Omega)}\|\varphi\|_p^p \nonumber\\
 &\leq \left[  \delta  + \frac{ C(N,p,r)S}{\delta^{N/(pr-N)}}K^{\frac{pr}{pr-N}}\right] \|\varphi\|_{1,p}^p = C_1,
\end{align}
where $C_1 > 0$ is independent of $V$.\\
Next, fix a $V_0\in \mrp$ such that $(V_0, m) \in \mathcal{W}_C $. Then we claim that
\begin{align}\label{max-8}
  T_0 =  \inf \left\{ \int_{\Omega}m|\varphi|^p \d x : (V_0,m) \in \mathcal{W}_C \right\} > 0.
\end{align}
 If not, there exists a sequence $\{(V_0,m_n)\} \in \mathcal{W}_C$ such that $\int_{\Omega}m_n|\varphi|^p \d x \to 0$ as $n\to \infty$. Since $\mathcal{W}_C$ is sequentially weakly compact (see \cite[Theorem 2.8.6]{Robert98}), hence there exists $(V_0,m_*) \in \mathcal{W}_C$ such that up to a subsequence $\{(V_0,m_n)\}$ converges weakly to $(V_0,m_*)$ in $ \mrp \times \mrp$. This implies that $m_n \rightharpoonup m_*$ in $\mrp$, and as $\varphi \in C_c^\infty(\Omega)$, we have 
 $$\int_{\Omega}m_n|\varphi|^p \d x\to \int_{\Omega}m_*|\varphi|^p \d x \text{ as } n \to \infty.$$ 
 Thus, the uniqueness of the limit implies that $\int_{\Omega}m_*|\varphi|^p \d x =0$, which is a contradiction as $m_* > 0$ a.e. in $\Omega$. Since the infimum as defined in \eqref{max-8} does not depend on $V_0$, hence we have \begin{align}\label{max-9}
    T_0=\inf \left\{ \int_{\Omega}m|\varphi|^p \d x : (V,m) \in \mathcal{W}_C \right\} > 0.
\end{align}
Therefore, by the definition of $\lambda_1(V, m)$ together with \eqref{max-7} and \eqref{max-8}, we get
\begin{align*}
        \lambda_1(V,m) &\leq \frac{\int_{\Omega} |\nabla \varphi|^p_A \d x + \int_{\Omega} V|\varphi|^p \d x}{\int_{\Omega} m|\varphi|^p} \nonumber\\
        &\leq \frac{ \Theta^{-p}_\Omega\int_{\Omega} |\nabla \varphi|^p \d x + \int_{\Omega} |V||\varphi|^p \d x}{\int_{\Omega} m|\varphi|^p} 
        %&\leq \frac{ S \Theta^{-p}_\Omega + \int_{\Omega} |V||\varphi|^p \d x}{\int_{\Omega} m|\varphi|^p} \nonumber\\
        \leq \frac{\Theta^{-p}_\Omega \|\varphi\|_{1,p}^p  + C_1}{T_0}. 
    \end{align*}
    Hence, the Lemma follows.
\end{proof}
As a consequence, we establish the existence of a maximizer for the functional $\lambda_1(V,m)$ over a weakly compact set $\mathcal{W}_C$.
\begin{thm}
 The functional $\lambda_1(V,m)$ admits a maximizer $(V^*,m^*) \in \mathcal{W}_C$. Furthermore, if $\mathcal{W}_C$ is strictly convex, then the maximizer is unique. 
\end{thm}
\begin{proof}
    Define \begin{align}\label{max-00}
        T_2 := \sup\{ \lambda_1(V,m) : (V,m) \in \mathcal{W}_{C}\}.
    \end{align}
    This is well-defined due to the Lemma \ref{bdd_lemma}. Then, there exists a sequence $\{(V_n,m_n)\}$ in $\mathcal{W}_{C}$ such that $\lambda_1(V_n,m_n) \to T_2$ as $n \to \infty$. Since $\mathcal{W}_C$ is sequentially weakly compact, hence 
 there exists $(V^*,m^*) \in \mathcal{W}_C$ such that up to a subsequence  $(V_n,m_n)$ converges weakly to $(V^*,m^*)$ in $\mrp \times \mrp$. This further implies that $V_n \rightharpoonup V^*$ and $m_n \rightharpoonup m^*$ in $\mrp$. On the other hand, for every $\varepsilon >0$, there exists $ u_\varepsilon \in \Wkp \setminus \{0\}$ such that
\begin{align}\label{max-4}
    \lambda_1(V^*,m^*) + \varepsilon  \geq \frac{ \int_{\Omega} |\nabla u_\varepsilon|^p_A \d x + \int_{\Omega}V^*|u_\varepsilon|^p \d x }{\int_{\Omega}m^*|u_\varepsilon|^p \d x}.
    % &= \lim_{n \to \infty} \frac{ \int_{\Omega} |\nabla w_k|^p_A \d x + \int_{\Omega}V_n|w_k|^p \d x }{\int_{\Omega}m_n|w_k|^p \d x}\nonumber\\
    % & \geq \lim_{n \to \infty} \lambda_1(V_n,m_n) = T_2.
\end{align}
Further, applying the Lemma \ref{convergence-2}, we have 
\begin{align}\label{max-2}
    \int_{\Omega} V_n|u_\varepsilon|^p \d x \to \int_{\Omega}V^*|u_\varepsilon|^p \d x \text{ as } n \to \infty,
\end{align}
and 
\begin{align}\label{max-3}
    \int_{\Omega} m_n|u_\varepsilon|^p \d x \to \int_{\Omega}m^*|u_\varepsilon|^p \d x \text{ as } n \to \infty.
\end{align}
Note that $\int_{\Omega}m^*|u_\varepsilon|^p \d x > 0$. Then, using \eqref{max-2} and \eqref{max-3} in \eqref{max-4}, we arrive at
\begin{align*}
     \lambda_1(V^*,m^*) + \varepsilon &\geq \lim_{n \to \infty} \frac{ \int_{\Omega} |\nabla u_\varepsilon|^p_A \d x + \int_{\Omega}V_n|u_\varepsilon|^p \d x }{\int_{\Omega}m_n|u_\varepsilon|^p \d x}
      \geq \lim_{n \to \infty} \lambda_1(V_n,m_n) = T_2.
\end{align*}
Thus, we get $\lambda_1(V^*,m^*) \geq T_2$. Also by the definition of supremum $\lambda_1(V^*,m^*) \leq T_2$. Therefore, we obtain $\lambda_1(V^*,m^*)=\sup\{ \lambda_1(V,m) : (V,m) \in \mathcal{W}_{C}\}$. 
Next, we show that the maximizer is unique under the hypothesis that $ \mathcal{W}_C$ is strictly convex subset. Let us suppose that $(V_1,m_1) \neq (V_2,m_2)$ are two distinct maximizers of $\lambda_1(V,m)$ in $\mathcal{W}_C$, i.e., 
\begin{align}\label{unique-00}
    \lambda_1(V_1,m_1) = \lambda_1(V_2,m_2) = \max\{\lambda_1(V,m) : (V,m) \in \mathcal{W}_C \} = T_2.
\end{align}
Let  $t \in (0,1)$. Since $\mathcal{W}_C$ is strictly convex subset, hence $(V_t,m_t) =  (1-t)(V_1,m_1) + t(V_2,m_2) \in Int (\mathcal{W}_C)$. This implies that $V_t = (1-t)V_1+tV_2$, and $m_t = (1-t)m_1+ tm_2$. Assume that $u_t$ is a positive  eigenfunction corresponding to the eigenvalue $\lambda_1(V_t,m_t)$ such that $\int_{\Omega}m_t|u_t|^p \d x = 1$. Then, we have 
\begin{align*}
    \lambda_1(V_t,m_t) = \int_{\Omega}|\nabla u_t|^p_A \d x + \int_{\Omega} V_t |u_t|^p \d x.
\end{align*}
Further, by using $ u_t$ as an admissible function in the definition of  $\lambda_1(V_1,m_1)$, we get
\begin{align}\label{unique-01}
   T_2 \int_{\Omega}m_1 |u_t|^p \d x = \lambda_1(V_1,m_1) \int_{\Omega}m_1 |u_t|^p \d x \leq \int_{\Omega}|\nabla u_t|^p_A \d x + \int_{\Omega} V_1 |u_t|^p \d x.
\end{align}
Parallelly, using $ u_t$ as an admissible function in the definition of $\lambda_1(V_2,m_2)$, we obtain
\begin{align}\label{unique-02}
   T_2 \int_{\Omega}m_2 |u_t|^p \d x = \lambda_1(V_2,m_2) \int_{\Omega}m_2 |u_t|^p \d x \leq \int_{\Omega}|\nabla u_t|^p_A \d x + \int_{\Omega} V_2 |u_t|^p \d x.
\end{align}
Multiplying \eqref{unique-01} by $(1-t)$, and \eqref{unique-02} by $t$, and then adding the resulting equations, we arrive at
\begin{align*}
    &T_2 \int_{\Omega} ((1-t)m_1+ tm_2)|u_t|^p \d x \leq \int_{\Omega}|\nabla u_t|^p_A \d x + \int_{\Omega} ((1-t)V_1+tV_2) |u_t|^p \d x  \\
    \implies &T_2 \int_{\Omega} m_t|u_t|^p \d x \leq \int_{\Omega}|\nabla u_t|^p_A \d x + \int_{\Omega} V_t |u_t|^p \d x \\
    \implies & T_2 \leq \lambda_1 (V_t,m_t).
\end{align*}
Therefore, by using \eqref{unique-00}, we obtain
\begin{align}\label{unique-03}
    \lambda_1(V_t,m_t) = T_2.
\end{align}
Let $u_1$ and $u_2$ be two positive eigenfunctions associated with eigenvalues $\lambda_1(V_1,m_1)$ and $\lambda_1(V_2,m_2)$, respectively, such that $\int_{\Omega}m_1|u_1|^p \d x = 1$ and $\int_{\Omega}m_2|u_2|^p \d x = 1$. We next claim that $u_1 = u_2 = u_t$. Suppose, for the sake of contradiction, that $u_t \neq u_1$ and $u_t \neq u_2$. Then from \eqref{unique-01} and \eqref{unique-02}, it follows that 
\begin{align}\label{unique-04}
    \int_{\Omega}|\nabla u_t|^p_A \d x + \int_{\Omega} V_1|u_t|^p \d x >  T_2  \int_{\Omega}m_1 |u_t|^p \d x,
\end{align}
and 
\begin{align}\label{unique-05}
    \int_{\Omega}|\nabla u_t|^p_A \d x + \int_{\Omega} V_2|u_t|^p \d x >   T_2  \int_{\Omega}m_2 |u_t|^p \d x.
\end{align}
Then by the virtue of inequalities \eqref{unique-04} and \eqref{unique-05}, it follows that
\begin{align*}
    \lambda_1(V_t,m_t) &= \int_{\Omega}|\nabla u_t|^p_A \d x + \int_{\Omega} V_t|u_t|^p \d x\\
    & =(1-t) \left[ \int_{\Omega}|\nabla u_t|^p_A \d x + \int_{\Omega} V_1|u_t|^p \d x\right] + t \left[ \int_{\Omega}|\nabla u_t|^p_A \d x + \int_{\Omega} V_2|u_t|^p \d x\right]\\
    & > (1-t) T_2\int_{\Omega}m_1 |u_t|^p \d x + t T_2\int_{\Omega}m_2 |u_t|^p \d x \\
    & = T_2 \int_{\Omega}m_t |u_t|^p \d x = T_2,
\end{align*}
which is a contradiction to \eqref{unique-03}. Therefore, we have $u_1 = u_2 = u_t$  in $\Wkp$. 

On the other hand, by the definition of weak solution \eqref{Weak}, we have 
\begin{align}\label{unique-06}
    \int_{\Omega}|\nabla u_1|^{p-2}_A\nabla u_1 \cdot \nabla \varphi \d x + \int_{\Omega} V_1|u_1|^{p-2}u_1 \varphi \d x = \lambda_1(V_1,m_1) \int_{\Omega}m_1|u_1|^{p-2}u_1 \varphi \d x,
\end{align}
and 
\begin{align}\label{unique-07}
     \int_{\Omega}|\nabla u_2|^{p-2}_A\nabla u_2 \cdot \nabla \varphi + \int_{\Omega} V_2|u_2|^{p-2}u_2 \varphi \d x = \lambda_1(V_2,m_2) \int_{\Omega}m_2|u_2|^{p-2}u_2 \varphi \d x,
\end{align} for all $ \varphi \in \Wkp$.
Subtracting  \eqref{unique-07} from \eqref{unique-06} and using the equalities $u_1 = u_2$,  and $\lambda_1(V_1,m_1) = \lambda_1(V_2,m_2)=T_2$, we find that
\begin{align*}
    \int_{\Omega}\left[(V_1 - V_2)-T_2(m_1-m_2)\right]|u_1|^{p-2}u_1 \varphi \d x =  0, \text{ for all } \varphi \in \Wkp.
\end{align*}
Since $u_1 > 0$, hence we have
\begin{align}\label{unique-08}
    V_1-V_2=T_2 (m_1-m_2) \text{ a.e. in } \Omega.
\end{align}
This implies that $m_1=m_2$ a.e. in $\Omega$ if and only if $V_1=V_2$ a.e. in $\Omega$. Since $(V_1,m_1) \neq (V_2,m_2)$, hence we have either (i) $V_1 = V_2$ and $m_1 \neq m_2$, or (ii) $V_1 \neq V_2$ and $m_1 = m_2$, or (iii) $V_1 \neq V_2$ and $m_1 \neq m_2$. Therefore, we must have $V_1 \neq V_2$ and $m_1 \neq m_2$. 

\noindent Since $(V_t,m_t) \in Int(\mathcal{W}_C)$, hence for a sufficiently small $\varepsilon>0$, $(V_t+\varepsilon,m_t) \in \mathcal{W}_C$. Let $u_\varepsilon$ be a positive eigenfunction associated with $\lambda_1(V_t+\varepsilon,m_t)$ such that $\int_{\Omega}m_t|u_\varepsilon|^p \d x = 1$. Then, we have
\begin{align}\label{unique-10}
     \lambda_1(V_t+\varepsilon,m_t)  %&= \int_{\Omega}|\nabla u_\varepsilon|^p_A \d x + \int_{\Omega}(V_t+\varepsilon)|u_\varepsilon|^p \d  x \nonumber\\
      = \int_{\Omega}|\nabla u_\varepsilon|^p_A \d x + \int_{\Omega}V_t|u_\varepsilon|^p \d  x + \varepsilon \int_{\Omega} |u_\varepsilon|^p \d  x. 
     % & \geq \lambda_1(V_t,m_t) +  \varepsilon \int_{\Omega} |u_\varepsilon|^p \d  x  > \lambda_1(V_t,m_t) = T_2,
 \end{align}
Moreover, using  $u_\varepsilon$ as an admissible in the definition of $\lambda_1(V_t,m_t)$, we get
\begin{align}\label{unique-09}
\lambda_1(V_t,m_t) & \leq     \frac{\int_{\Omega}|\nabla u_\varepsilon|^p_A \d x + \int_{\Omega}V_t|u_\varepsilon|^p \d x}{\int_{\Omega}m_t|u_\varepsilon|^p \d x }  =  \int_{\Omega}|\nabla u_\varepsilon|^p_A \d x + \int_{\Omega}V_t|u_\varepsilon|^p \d x.
\end{align}
Since $u_\varepsilon$ is a positive eigenfunction, hence using the inequality \eqref{unique-09} in \eqref{unique-10}, we arrive at
\begin{align*}
         \lambda_1(V_t+\varepsilon,m_t)  \geq \lambda_1(V_t,m_t) +  \varepsilon \int_{\Omega} |u_\varepsilon|^p \d  x  > \lambda_1(V_t,m_t) = T_2,
\end{align*}
which contradicts \eqref{max-00}. Therefore, we must have $(V_1,m_1) = (V_2,m_2)$, which establishes the uniqueness of maximizer.
\end{proof}
We now discuss the existence of a minimizer of the functional $\lambda_1(V,m)$ over the weakly compact set $\mathcal{W}_C$. To establish this result, we derive the following lower bound of the functional $\lambda_1(V,m)$ over $\mathcal{W}_C$.
\begin{lem}\label{bdd-below}
There exists $C(N,p,r,S,K)> 0$ such that $$\lambda_1(V, m) \geq  -C, $$
for all $(V, m)\in \mathcal{W}_C$.
\end{lem}
\begin{proof}
Let $(V, m) \in \mathcal{W}_C$. Assume that $u$ is an eigenfunction associated with $\lambda_1(V,m)$ such that $\|u\|_{1,p} = 1$. Then 
\begin{align*}
    \lambda_1(V,m) = \frac{\int_{\Omega} |\nabla u |^p_A \d x + \int_{\Omega} V |u|^p \d x }{\int_{\Omega}m|u|^p \d x }.
\end{align*}
    By the Morrey–Adams Theorem \ref{M-A-THM} (for a fix $\delta > 0$), the Poincar\'e inequality, and the bound $\Vert{}V\Vert{}_\mrp \leq K$, together with the definition of $T_0$ in \eqref{max-9}, we get
 \begin{align*}
    \lambda_1(V,m) &= \frac{\int_{\Omega} |\nabla u |^p_A \d x + \int_{\Omega} V|u|^p \d x }{\int_{\Omega}m|u|^p \d x } \geq  \frac{  -\int_{\Omega} |V||u|^p \d x }{\int_{\Omega}m|u|^p \d x } \geq -\frac{1}{T_0}\int_{\Omega} |V||u|^p \d x\\
    & \geq \frac{1}{T_0} \left[ - \delta \|u\|_{1,p}^p -  \frac{C(N,p,r)}{\delta^{N/(pr-N)}} \|V\|^{pr/(pr-N)}_{M^r(p;\Omega)} \|u\|_p^p\right]\\
    & \geq \frac{1}{T_0}\left[ - \delta -  \frac{C(N,p,r)}{\delta^{N/(pr-N)}} K^{pr/(pr-N)} S\right]  = -\frac{C_1 (N, p, r, S, K)}{T_0}.
\end{align*}
    Hence the Lemma follows. 
\end{proof}
As a result, we have: 
\begin{thm}\label{mini-thm}
The functional $\lambda_1(V,m)$ admits a  minimizer $(V_*,m_*) \in \mathcal{W}_C$.
\end{thm}
\begin{proof}
    Let $T_1 := \inf \{ \lambda_1(V,m) : (V, m) \in \mathcal{W}_{C}\}$, which is well-defined by the Lemma \ref{bdd-below}. By the definition of infimum, there exists a sequence $\{ (V_n,m_n)\}$ in $\mathcal{W}_C$ such that $\lambda_1(V_n,m_n) \to T_1$ as $n \to \infty$.  Let  $u_n$ be the eigenfunction corresponding to the  eigenvalue $\lambda_1(V_n,m_n)$ such that $ \int_{\Omega} m_n|u_n|^p \d x = 1$. Then 
        $$  \int_{\Omega} |\nabla u_n|^p_A \d x + \int_{\Omega} V_n|u_n|^p \d x =\lambda_1(V_n,m_n). $$
This implies that
$$\int_{\Omega} |\nabla u_n|^p_A \d x + \int_{\Omega} V_n|u_n|^p \d x \to T_1 \text{ as } n \to \infty.$$
    Since $\mathcal{W}_C$ is sequentially weakly compact, hence 
 there exists $(V_*,m_*) \in \mathcal{W}_C$ such that up to a subsequence  $(V_n,m_n) $ converges weakly to $(V_*,m_*)$ in $\mrp \times \mrp$. Hence, we have  $V_n \rightharpoonup V_*$ and $m_n \rightharpoonup m_*$ in $\mrp$.  Since $\{\lambda_1(V_n,m_n)\}$ is bounded (i.e., $|\lambda_1(V_n,m_n) |\leq C $), thus by the Lemma \ref{bdd-lema}, we get 
\begin{align*}
    \|u_n\|_{1,p} \leq C_1 \lambda_1(V_n,m_n) + C_2 \int_{\Omega} m_n|u_n|^p \d x \leq C_1 C + C_2.
\end{align*}
This shows that $\{u_n\}$ is bounded in $\Wkp$. So, there exists  $u_* \in \Wkp$ such that $u_n \rightharpoonup u_*$ in $\Wkp$ (up to a subsequence).  
Thus, by using the Lemma  \ref{convergence-3}, we have
\begin{align}\label{mini-1}
    \int_{\Omega} V_n|u_n|^p \d x \to \int_{\Omega}V_*|u_*|^p \d x \text{ as } n \to \infty,
\end{align}
and 
\begin{align}\label{mini-2}
    \int_{\Omega} m_n|u_n|^p \d x \to \int_{\Omega}m_*|u_*|^p \d x \text{ as } n \to \infty.
\end{align}
Since $\int_{\Omega}m_n|u_n|^p \ dx = 1$, hence we get $\int_{\Omega} m_*|u_*|^p \d x = 1$. Now,  choosing $u_*$ as an admissible function in the definition of  $\lambda_1(V_*,m_*)$ and then using the  weak lower semicontinuity of $\int_{\Omega} |\nabla u|^p_A \d x$ together with \eqref{mini-1}, we get
    \begin{align*}
        \lambda_1(V_*,m_*) \leq \int_{\Omega} |\nabla u_*|^p_A \d x + \int_{\Omega} V_*|u_*|^p \d x &\leq \liminf_{n \to \infty} \left( \int_{\Omega} |\nabla u_n|^p_A \d x + \int_{\Omega} V_n|u_n|^p \d x \right) = T_1.
    \end{align*}
This concludes the proof.
% Thus, $V_*$ and $m_*$ is a minimizer of $\lambda_1(V,m)$ in $\mathcal{W}_{C}$ and $u_*$ is an eigenfunction associated with $\lambda_1(V_*,m_*)$ for the problem \eqref{EVP} with potential $V_*$ and weight $ m_*$. 
\end{proof}

\begin{rk} It is worth noting that the minimizer of $\lambda_1(V,m)$ over $\mathcal{W}_C$ is not necessarily unique. Reader can refer to Biswas et. al \cite[Remark 7]{Biswas23} for the case $A = I_N$ and $V, m$ belongs to $L^{N/p}(\Omega)$. We also mention the work of Pezzo and Bonder \cite{Pezzo10} for the case $A= I_N$ and $V,m$ lies in $L^q(\Omega)$ with $q > \frac{N}{p}$, where they proved the existence of a minimizing pair $(V_*,m_*)$ over a rearrangement class of fixed $(V_0, m_0)$ such that $V_*$  and $m_*$ are the unique minimizer of  the functionals $ V \mapsto \int_{\Omega}V|u_*|^p \d x $ and $ m \mapsto \int_{\Omega}m|u_*|^p \d x $, respectively. 
\end{rk}
To the best of our knowledge, no explicit examples demonstrating the non-uniqueness of minimizer for the functional $\lambda_1(V,m)$ exist in the literature. Thus, we present the following example.
\begin{exmp} 
Let $\Omega = B_1(0) \subset \mathbb{R}^2$ be the unit disk,  $ 1 < p < 2$. For $ a \in [-1,1]$, define
$$V_a(x,y) := \begin{cases}
 a\,; & \text{ when } x>0,\\
 -a\,; & \text { when } x<0.
\end{cases}$$
Then $V_a(x,y) \in L^\infty (\Omega) \subset \mrp$ for every $a \in [-1,1]$. Consider the weight function $m:=1$ on $\Omega$. Take $\mathcal{W}_C := \{ (V_a,1) : a \in [-1,1] \}$. Let $T : [-1,1] \to \mrp \times \mrp $ be the map defined by  $ T(a) = (V_a,1) $. Notice that 
\begin{align*}
     \|T(a) - T(b)\|_{\mrp \times \mrp}  %= \|(V_a,1)-(V_b,1)\|_{\mrp \times \mrp}\\
     &= \|(V_a-V_b,0)\|_{\mrp \times \mrp} = \|(V_{a-b},0)\|_{\mrp \times \mrp}\\
     &= \|V_{a-b}\|_{\mrp} = \sup_{\substack{z \in \Omega \\ 0 < s < \operatorname{diam}(\Omega)}} \frac{1}{s^{2/r'}} \int_{\Omega \cap B_s(z)} |V_{a-b}(z)| \d z \\
     & \leq \sup_{\substack{z \in \Omega \\ 0 < s < \operatorname{diam}(\Omega)}} \frac{1}{s^{2/r'}}  |\Omega \cap B_s(z)| |a-b|\\
    % &\leq \sup_{\substack{z \in \Omega \\ 0 < s < \operatorname{diam}(\Omega)}} \frac{1}{s^{2/r'}}  \pi s^2 |a-b|\\
     & \leq \sup_{\substack{z \in \Omega \\ 0 < s < \operatorname{diam}(\Omega)}} \pi s^{2/r} |a-b|  \leq  2^{2/r} \pi |a-b|.
\end{align*}
This shows that the map $T$ is continuous and, hence $\mathcal{W}_C = Im(T)$ is a compact subset of $\{ (V, m)\in \mrp \times \mrp : m > 0 \text{ a.e.  in } \Omega \}$. Now by Theorem \ref{mini-thm}, there exists $a^*\in [-1, 1]$ such that  $(V_{a^*},1)$ is a minimizer of $\inf \{ \lambda_1(V,1) : (V, 1) \in \mathcal{W}_{C}\}$. 

We now claim that $a^*\neq 0$. To prove this, we show that $\lambda_1(V_a,1) < \lambda_1(V_0,1)=\lambda_1(0,1)$ for every $a \neq 0$. Without loss of generality, let us assume that $a >0$. Let $u_0$ be a positive eigenfunction associated with $\lambda_1(V_0,1)$ such that $\|u_0\|_p = 1$. Denote $\Omega_R := \{ (x,y) \in \Omega : x > 0\}$ and $\Omega_L := \{ (x,y) \in \Omega : x < 0 \}$. Let $\varphi \in C_c^\infty(\Omega)$ such that $\varphi(-x,y) = -\varphi(x,y)$ with  $ \varphi > 0$ in $\Omega_L$, and $ \varphi < 0 $  on $\Omega_R$. For a small enough $\varepsilon > 0$, define
$$u_\varepsilon(x,y) = u_0(x,y) + \varepsilon \varphi(x,y). $$
Then, using the Taylor expansion, we get 
$$|u_\varepsilon |^p = |u_0 + \varepsilon \varphi |^p = u_0^p + \varepsilon p  u_0^{p-1}\varphi + o(\varepsilon^2).$$
Denote $f(x,y):=  u_0^{p-1}(x,y)\varphi(x,y)$. Then $f(-x,y) = -f(x,y)$ and $\iint_{\Omega_L} f(x,y) \d x\d y=\iint_{\Omega_R} f(-x,y) \d x \d y =  -\iint_{\Omega_R} f(x, y)\d x \d y $. Thus $\iint_\Omega f(x,y)\d x \d y= \iint_{\Omega_L} f(x,y)\d x \d y + \iint_{\Omega_R} f(x,y)\d x \d y =0$. Using this, we obtain
\begin{align}\label{examp-8}
    \|u_\varepsilon\|_p^p  = 1 +  o(\varepsilon^2).
\end{align}
Since $V_a(-x,y)=-V_a(x, y)$, similarly we have
\begin{align*}
\iint_{\Omega} V_a |u_\varepsilon|^p  \d x \d y&= \iint_{\Omega}V_a u_0^p \d x \d y + \varepsilon p \iint_{\Omega} V_a u_0^{p-1} \varphi \d x \d y  + o(\varepsilon^2)\\
& =   \varepsilon p \iint_{\Omega} V_a u_0^{p-1} \varphi \d x \d y  + o(\varepsilon^2).
\end{align*} 
Now, 
\begin{align*}
   \varepsilon p\iint_{\Omega} V_au_0^{p-1} \varphi \d x \d y &= \varepsilon p\iint_{\Omega_R} a u_0^{p-1} \varphi \d x \d y + \varepsilon p \iint_{\Omega_L} (-a)u_0^{p-1} \varphi \d x \d y\\
   &= -2 p a \varepsilon \iint_{\Omega_L}  u_0^{p-1} \varphi \d x \d y = - C_1 \varepsilon, 
\end{align*}
where $C_1 =  2 p a \iint_{\Omega_L}  u_0^{p-1} \varphi \d x \d y$ > 0.
Thus, we get 
\begin{align}\label{examp-7}
    \iint_{\Omega} V_a |u_\varepsilon|^p  \d x \d y = - C_1 \varepsilon +  o(\varepsilon^2).
\end{align}
Furthermore, we have $$| \nabla u_0 + \varepsilon \nabla \varphi |^p = |\nabla u_0|^p +  \varepsilon p|\nabla u_0|^{p-2} \nabla u_0 \cdot \nabla \varphi +  o(\varepsilon^2), $$
and consequently, 
\begin{equation}\label{exx1}
\iint_{\Omega} |\nabla u_\varepsilon |^p \d x \d y =\iint_{\Omega} |\nabla u_0 |^p \d x \d y + p \varepsilon \iint_{\Omega} |\nabla u_0|^{p-2} \nabla u_0 \cdot \nabla \varphi \d x \d y + o(\varepsilon^2).
\end{equation}
Note that $\nabla u_0(x,y) = u'_0(\tilde{r}) (\frac{x}{\tilde{r}},\frac{y}{\tilde{r}})$, where $\tilde{r} = \sqrt{x^2+y^2}$. Then 
\begin{align*}
    \frac{\partial u_0}{\partial x}(-x,y)=u'_0(\tilde{r})\left(\frac{-x}{\tilde{r}}\right) = -u'_0(\tilde{r})\left(\frac{x}{\tilde{r}}\right) = -\frac{\partial u_0}{\partial x}(x,y),
\end{align*}
and similarly, we get 
\begin{align*}
    \frac{\partial u_0}{\partial y}(-x,y)=u'_0(\tilde{r})\left(\frac{y}{\tilde{r}}\right) =  \frac{\partial u_0}{\partial y}(x,y).
\end{align*}
On the other hand,  since $\varphi(-x,y) = - \varphi(x,y)$, hence
\begin{align}\label{examp-1}
    \frac{\partial }{\partial x}(\varphi(-x,y)) = -\frac{\partial \varphi}{\partial x}(x,y).
\end{align}
 By the chain rule,
\begin{align}\label{examp-2}
    \frac{\partial }{\partial x}(\varphi(-x,y)) =  -\frac{\partial \varphi }{\partial x}(-x,y).
\end{align}
Thus from \eqref{examp-1} and \eqref{examp-2}, we have 
\begin{align*}
    \frac{\partial \varphi}{\partial x} (-x,y) =  \frac{\partial \varphi}{\partial x} (x,y).
\end{align*}
Similarly, we get 
\begin{align*}
    \frac{\partial \varphi}{\partial y} (-x,y) =  -\frac{\partial \varphi}{\partial y} (x,y).
\end{align*}
Thus, the integrand $ |\nabla u_0|^{p-2} \nabla u_0 \cdot \nabla \varphi $ is an odd function in $x$. Therefore, 
\begin{align*}
    \iint_{\Omega} |\nabla u_0|^{p-2} \nabla u_0 \cdot \nabla \varphi  \d x \d y = 0.
\end{align*}
Hence, from \eqref{exx1}, we arrive at
\begin{align}\label{examp-10}
    \iint_{\Omega} |\nabla u_\varepsilon |^p \d x \d y = \iint_{\Omega} |\nabla u_0|^p \d x \d y +  o(\varepsilon^2).
\end{align}
Next, using \eqref{examp-8}, \eqref{examp-7} and \eqref{examp-10} in the  definition of $\lambda_1(V,m)$ for $A = I_2$ together with  potential $V = V_a$ and weight $m = 1$, we deduce that
\begin{align}\label{examp-11}
    \lambda_1(V_a,1) \leq \frac{\iint_{\Omega} (|\nabla u_\varepsilon |^p + V_a |u_\varepsilon|^p ) \d x \d y}{\|u_\varepsilon\|^p_p} &= \frac{\lambda_1(0,1) - C_1 \varepsilon  + o(\varepsilon^2)}{1 + o(\varepsilon^2)} \nonumber\\
    &\leq \lambda_1(0,1) - C_1 \varepsilon  + o(\varepsilon^2).
\end{align}
Choose $\varepsilon > 0$ sufficiently small enough so that $ C_1 \varepsilon > o(\varepsilon^2)$. Then from  \eqref{examp-11}, we find that
$$\lambda_1(V_a,1) \leq \lambda_1(0,1) - C_1 \varepsilon  + o(\varepsilon^2) < \lambda_1(0,1).  $$

We finally show that $(V_{-a^*},1)$ is also a minimizer of $\lambda_1(V,m)$ on $\mathcal{W}_C$. For that, consider the reflection operator $R : \Wkp \to \Wkp$ defined by $R(u(x,y)) = u(-x,y)$. Then, a change of variable implies
\begin{align}\label{examp-13}
    \iint_{\Omega} |R(u(x,y))|^p \d x  \d y =  \iint_{\Omega} |u(-x,y)|^p \d x \d y = \iint_{\Omega} |u(x,y)|^p \d x  \d y,
\end{align}
and 
 \begin{align}\label{examp-14}
     \iint_{\Omega} V_{-a^*}(x,y) |u(x,y)|^p \d x \d y &= \iint_{\Omega} V_{a^*}(-x,y) |u(x,y)|^p \d x \d y \nonumber\\
    % & = \iint_{\Omega} V_{a^*}(x,y) |u(-x,y)|^p \d x \d y = \iint_{\Omega} V_{a^*}(x,y) |R(u(x,y))|^p \d x \d y\nonumber\\
     & = \iint_{\Omega} V_{a^*}(x,y) |R(u(x,y))|^p \d x \d y.
 \end{align}
Moreover, \begin{align}\label{examp-15}
    \iint_{\Omega} |\nabla Ru(x,y)|^p \d x \d y &=  \iint_{\Omega} |\nabla u(-x,y)|^p \d x \d y= \iint_{\Omega} |\nabla u(x,y) |^p \d x \d y.
    % & = \iint_{\Omega} |\nabla u(z,w) |^p \d z \d w \nonumber\\
    % & = \iint_{\Omega} |\nabla u(x,y) |^p \d x \d y
\end{align}
Substituting \eqref{examp-13}, \eqref{examp-14} and \eqref{examp-15}  into the definition of $\lambda_1(V_{-a^*}, 1)$ with $A = I_2$, and using the bijectivity of reflection operator $R$, we obtain
\begin{align*}
    \lambda_1(V_{-a^*},1) &= \inf_{u \in \Wkp\setminus \{0\}} \frac{\iint_{\Omega}|\nabla u(x,y)|^p \d x \d y + \iint_{\Omega}V_{-a^*}(x,y) |u(x,y)|^p \d x \d y}{\iint_{\Omega}|u(x,y)|^p \d x \d y} \\
    & = \inf_{u \in \Wkp\setminus \{0\}} \frac{\iint_{\Omega}|\nabla R(u(x,y))|^p \d x \d y + \iint_{\Omega}V_{a^*}(x,y) |R(u(x,y))|^p \d x \d y}{\iint_{\Omega}|R(u(x,y))|^p \d x \d y} \\
    & = \lambda_1(V_{a^*},1).
\end{align*}
Thus, minimizers of $\{ \lambda_1(V,1) : (V, 1) \in \mathcal{W}_{C}\}$ is not unique.
\end{exmp}

\section{Derivative of the functional $\lambda_1(V,m)$}
In this section, we address the differentiability of the principal eigenvalue $\lambda_1(V, m)$ with respect to regular perturbations of potential $V$ and weight $m$ within a closed subspace  $M^r_0(p;\Omega) := \overline{C_c^\infty(\Omega)}$ of the Morrey space $\mrp$ where $r>\frac{N}{p}$. Our finding build upon and extend the results of Pezzo and Bonder \cite{Pezzo10}. 

Throughout this section, we assume potential $V$, and positive weight $m$ lies in closed subspace $M^r_0(p; \Omega)$. Let $W$ be a regular (smooth), globally Lipschitz vector field compactly supported in $\Omega$, with $ \|\nabla W\|_\infty < \infty$ and $\|W\|_\infty < \infty$.  Let the flow $\psi_t : \mathbb{R}^N \to \mathbb{R}^N$ be generated by the autonomous initial value problem
    \begin{align}\label{flow}
    \begin{split}
   \frac{\d}{\d t}\psi_t(x) &= W(\psi_t(x));\,\, t >0,\\
    \psi_0(x) &= x;\,\, x \in \mathbb{R}^N.
    \end{split}
\end{align}
Note that $\psi_t$ is a $C^1$-diffeomorphisms with $\psi^{-1}_{t}=\psi_{-t}$. Also, as $\operatorname{supp}(W) \subset \Omega$, hence we have $\psi_t(\Omega) = \Omega$, $\psi_t(x)=x$ for $x\in \mathbb{R}^N\setminus \Omega$, and $\psi_t^{-1}(\Omega) = \Omega$ for all  $t \geq 0$. Furthermore, by Taylor expansion, $\psi_t(x) = x + t W(x) + o(t^2) \text{ for all } x \in \mathbb{R}^N$.

For $V, m (>0) \in M^r_0(p; \Omega)$, let us define $V_t := V \circ \psi_t^{-1}$ and $m_t = m \circ \psi_t^{-1}$. Suppose that $u_t$ is a positive eigenfunction associated with   $\lambda_1(V_t,m_t)$ such that $\int_{\Omega}m_t|u_t|^p \d x = 1$. Then 
\begin{align}\label{eigg}
    \lambda_1(t):=\lambda_1(V_t,m_t) = \int_{\Omega}|\nabla u_t|^p \d x + \int_{\Omega} V_t |u_t|^p \d x.
\end{align}

\begin{lem}\label{convergence-m} 
 Let $ V \in M^r_0(p;\Omega) $, then $V_t \rightarrow V$ in  $M^r_0(p;\Omega)$ as $ t \to 0^+$.
\end{lem}
\begin{proof}
  Let  $\varepsilon > 0$. Then by the definition of $M^r_0(p;\Omega)$,  there exists $V_\varepsilon \in C_c^\infty(\Omega)$ such that
  \begin{align}\label{conv-gt-0}
      \|V_\varepsilon - V\|_{\mrp} < \varepsilon.
  \end{align}
  Denote $(V_\varepsilon)_t: = V_\varepsilon \circ \psi^{-1}_t $. Then,
  \begin{align}\label{conv-gt-1}
      \|V_t-V\|_{\mrp} \leq \|V_t - (V_\varepsilon)_t \|_{\mrp} + \|(V_\varepsilon)_t - V_\varepsilon \|_{\mrp} + \|V_\varepsilon- V \|_{\mrp}.
  \end{align}
Since $V_\varepsilon \in C_c^\infty(\Omega)$, and the inverse flow map $\psi_t^{-1}$ converges uniformly to the identity map on $\bar{\Omega}$ as $t \to 0^+$, it follows that $(V_\varepsilon)_t \to V_\varepsilon$ uniformly on $\Omega$ (see Pezzo and Bonder \cite{Pezzo10}). Thus, there exists $t_1 > 0$ such that 
$\| (V_\varepsilon)_t - V\|_\infty < \varepsilon $ for all $ 0 < t < t_1$.
Now,
  \begin{align}
      \|(V_\varepsilon)_t &- V_\varepsilon\|_{\mrp} = \sup_{\substack{y \in \Omega \\ 0 < s < \operatorname{diam}(\Omega)}} \frac{1}{s^{N/r'}} \int_{\Omega \cap B_s(x)} |(V_\varepsilon)_t - V_\varepsilon| \d y \nonumber\\
     % & = \sup_{\substack{y \in \Omega \\ 0 < s < \operatorname{diam}(\Omega)}} \frac{1}{s^{N/r'}} \|(V_\varepsilon)_t - V_\varepsilon\|_{\infty} |\Omega \cap B_s(x)|  \nonumber\\
      & \leq  \sup_{\substack{y \in \Omega \\ 0 < s < \operatorname{diam}(\Omega)}} C \frac{s^N}{s^{N/r'}} \|(V_\varepsilon)_t - V_\varepsilon\|_{\infty}   
     % &= \sup_{\substack{y \in \Omega \\ 0 < s < \operatorname{diam}(\Omega)}} C s^{N/r'} \|(V_\varepsilon)_t - V_\varepsilon\|_{\infty}    \nonumber\\
       \leq C ( \operatorname{daim} (\Omega) )^{N/r} \|(V_\varepsilon)_t - V_\varepsilon\|_{\infty}. \nonumber
  \end{align}
  Hence for all $ 0 < t < t_1$, we obtain
  \begin{align}\label{conv-gt-2}
      \|(V_\varepsilon)_t-V_\varepsilon\|_{\mrp} < C_1\varepsilon,
  \end{align}
  where $C_1 = C ( \operatorname{daim} (\Omega) )^{N/r} $  is a positive constant independent of $t$. Next, applying the change of variable  $ y = \psi_t(z)$, we have
\begin{align}\label{conv-gt-4}
    \int_{\Omega \cap B_s(x)} |V_t(y)| \d y =   \int_{\Omega \cap B_s(x)} |V(\psi_t^{-1}(y))| \d y = \int_{\psi_t^{-1}(\Omega \cap B_s(x))} |V(z)| |J_{\psi_t}| \d z,
 \end{align}
where $J_{\psi_t}$ is the Jacobian of $\psi_t$, given by $J_{\psi_t} = 1 + t \operatorname{div} W+  o(t^2)$. Since $W$ is a fix, smooth vector field with compact support in $\Omega$, hence $\operatorname{div} W$ is a continuous function on a compact subset  of $\Omega$. Thus $\displaystyle{ L_W:  =\sup_{x \in \Omega} \operatorname{div}  W < \infty}$, and consequently, there exists  $t_2 >0$ small enough such that $|J_{\psi_t}| \leq 1 + t L_W + o(t^2)$ for all $ 0 < t < t_2$. 
Using  $\psi_t^{-1}(\Omega) = \Omega$ and  $\psi_t^{-1}(B_s(x)) \subset B_{s_t}(\psi_t^{-1}(x))$, where $s_t = s(1 + tL_W+ o(t^2))$ in \eqref{conv-gt-4}, we arrive at
  \begin{align}\label{conv-gt-5}
       \int_{\Omega \cap B_s(x)} |V_t(y)| \d y \leq (1 + t L_W + o(t^2)) \int_{\Omega \cap B_{s_t}(\psi_t^{-1}(x))}| V(z)| \d z.
  \end{align}
Multiplying both sides by $s^{-N/r'}$ in \eqref{conv-gt-5}, we get
  \begin{align*}
     s^{-N/r'}  \int_{\Omega \cap B_s(x)} |V_t(y)| \d y &\leq (1 + t L_W + o(t^2)) s^{-N/r'} \int_{\Omega \cap B_{s_t}(\psi_t^{-1}(x))}| V(z)| \d z \\
     &=(1 + t L_W + o(t^2)) s_t^{-N/r'} \left(\frac{s_t}{s}\right)^{N/r'}  \int_{\Omega \cap B_{s_t}(\psi_t^{-1}(x))}| V(z)| \d z \\
     & = ( 1+ tL_W + o(t^2))^{1 + N/r'} s_t^{-N/r'}  \int_{\Omega \cap B_{s_t}(\psi_t^{-1}(x))}| V(z)| \d z.
  \end{align*}
Thus, we find that
\begin{align*}
    \|V_t\|_{\mrp}&=  \sup_{\substack{y \in \Omega \\ 0 < s < \operatorname{diam}(\Omega)}} \frac{1}{s^{N/r'}} \int_{\Omega \cap B_s(x)} |V_t(y)| \d y \nonumber\\
    & \leq \left( 1+ tL_W \left(1+\frac{N}{r'}\right) + o(t^2)\right) \|V\|_{\mrp}.
    \end{align*}
    Therefore, for $ 0 < t < t_2$ (where $t_2$ is small enough), we obtain
    \begin{align}\label{conv-gt-3}
         \|V_t\|_{\mrp} \leq C_W \|V\|_{\mrp},
    \end{align}
    where $C_W$ is a positive constant independent of $t$. Finally, using \eqref{conv-gt-0}, \eqref{conv-gt-2} and \eqref{conv-gt-3} in \eqref{conv-gt-1}, we have 
\begin{align*}
     \|V_t-V\|_{\mrp} \leq (1 + C_W + C_1) \varepsilon,
\end{align*}
for all $ 0 < t < \min\{t_1,t_2\}$. Therefore, the Lemma follows.
\end{proof}

As a consequence of the Lemma  \ref{convergence-m} and the Theorem \ref{continuous-thm}, we have:

\begin{cor}\label{continous-t}
    Let $V, m(>0) \in M^r_0(p;\Omega)$. Then $\lambda_1(t)$ is continuous at $t = 0$, i.e.,
    $$ \lambda_1(V_t,m_t) \to \lambda_1(V,m) \text{ as }  t \to 0^+.$$
\end{cor}

\begin{lem}\label{convergence-u}
    Let $V, m(>0) \in M^r_0(p;\Omega)$. Let $u_t$ be a  positive eigenfunction associated with $\lambda_1(V_t,m_t)$ with $t>0$ such that $ \int_{\Omega}m_t|u_t|^p \d x = 1$. Then  
    $$u_t \to u_0 \text{ in } \Wkp \text{ as } t\rightarrow 0^+, $$
    where $u_0$ is the positive eigenfunction associated with $\lambda_1(V, m)$ such that $\int_{\Omega}m|u_0|^p \d x = 1$. Furthermore, $u_t \circ \psi_t \rightarrow u_0$ in $\Wkp$ as $t \to 0^+$.
\end{lem}
\begin{proof}
   By the Corollary \ref{continous-t}, it follows that $\lambda_1(V_t,m_t) \leq C$, for some positive constant $C > 0$. Consequently, applying the Lemma \ref{bdd-lema}, we have $\{u_t\}$ is bounded in $\Wkp$. Therefore,  for a given $\{t_n\}$, the sequence  $\{u_{t_n}\}$ is bounded in $\Wkp$.
   Thus there exists  $u_0\in \Wkp$ such that up to a subsequence (still denoted by $u_{t_n}$) $ u_{t_n} \rightharpoonup u_0  \text{ in } \Wkp$.
Also, by the Lemma \ref{convergence-m}, $V_{t_n} \to V$ and $m_{t_n} \to m$ in $M^r_0(p;\Omega)$ as $t_n \to 0^+$. Then, using the Lemma \ref{convergence-3}, we have 
\begin{equation}\label{contii-1}
\left\{
\begin{aligned}
\int_{\Omega}V_{t_n}|u_{t_n}|^p \d x \to \int_{\Omega}V|u_0|^p \d x, \\
\int_{\Omega}m_{t_n}|u_{t_n}|^p \d x \to \int_{\Omega}m|u_0|^p \d x,
\end{aligned}
\right.
\end{equation}
as $t_n \to 0^+$. Since $\int_{\Omega}m_{t_n}|u_{t_n}|^p \d x = 1 $, hence we get $\int_{\Omega}m|u_0|^p \d x  = 1$. Then, using \eqref{contii-1} we have 
\begin{align*}
    \lambda_1(V,m) &= \lim_{ t_n \to 0^+} \lambda_1(V_{t_n},m_{t_n}) \nonumber \\
    &= \lim_{ t_n \to 0^+} \left( \int_{\Omega} |\nabla u_{t_n}|^p_A \d x + \int_{\Omega} V_{t_n}|u_{t_n}|^p \d x\right) \nonumber \\
    & \geq \liminf_{ t_n \to 0^+} \left( \int_{\Omega} |\nabla u_{t_n}|^p_A \d x + \int_{\Omega} V_{t_n}|u_{t_n}|^p \d x\right) \nonumber\\
    & \geq \int_{\Omega}|\nabla u_0|^p_A \d x + \int_{\Omega} V|u_0|^p \d x \geq \lambda_1(V,m). 
\end{align*}
This shows that $u_0$ is an eigenfunction corresponding to $\lambda_1(V,m)$, and as each $u_{t_n}$ is positive, $u_0$ is positive. Moreover, the convergence $\int_{\Omega} |\nabla u_{t_n}|^p_A \d x \to \int_{\Omega}|\nabla u_0|^p_A \d x $ as $ t_n \to 0^+$, which implies that $u_{t_n} \to u_0$  in $\Wkp$ as $ t_n \to 0^+$.

Since $u_t \to u_0$  in $\Wkp$ as $t\rightarrow 0^+$, and the flow map $\psi_t \to \text{Id}$ uniformly on $\bar{\Omega}$, hence we have  $u_t \circ \psi_t \to u_0$ in $\Wkp$ as $ t \to 0^+$.
\end{proof}
We thus obtain the following result regarding the differentiability of the functional $\lambda_1(t)$. 
\begin{thm} 
 Let $V, m(>0) \in M^r_0(p;\Omega)$. Consider the eigenfunction $u_0$ associated with $\lambda_1(V,m)$, as given in the Lemma \ref{convergence-u}. Then
$\lambda_1(t):=\lambda_1(V_t,m_t)$ is differentiable at $t = 0$ with derivative
\begin{align*}
    \frac{\d \lambda_1(t) }{\d t} \Bigg|_{t=0} &= \int_{\Omega} \left[|\nabla u_0|^p_A  + V|u_0|^p \right]\operatorname{div} W\d x - p \int_{\Omega}|\nabla u_0|^{p-2}_A \left<A \nabla u_0, W' \nabla u_0 \right> \d x  \\
    &\hspace{5cm}- \lambda_1(V,m) \int_{\Omega}m |u_0|^p  \operatorname{div} W \d x,
\end{align*}
where $W'$ denotes the Jacobian matrix of $W$.
\end{thm}
\begin{proof}
Define $g_t := u_0 \circ \psi_t^{-1}$.  Then 
\begin{align}\label{deriv-1}
    \lambda_1(V_t,m_t) \int_{\Omega}m_t|g_t|^p \d x \leq \int_{\Omega} |\nabla g_t|_A^p \d x + \int_{\Omega} V_t |g_t|^p \d x. 
\end{align}
Now, the change of variables yields
    \begin{align}\label{deriv-2}
        \int_{\Omega}m_t |g_t|^p \d x &= \int_{\Omega} m |u_0|^p |J \psi_t|  \d x \nonumber\\
        & = \int_{\Omega} m |u_0|^p ( 1 + t \operatorname{div} W  + o(t^2) )  \d x \nonumber\\
        & = 1 + t \int_{\Omega} m|u_0|^p \operatorname{div}W \d x + o(t^2).
    \end{align}
    Similarly, we obtain
    \begin{align}\label{deriv-3}
  \int_{\Omega}V_t |g_t|^p \d x & = \int_{\Omega} V |u_0|^p + t \int_{\Omega} V|u_0|^p \operatorname{div}W \d x + o(t^2).
    \end{align}
    Note that $\nabla u_0(x) =  [\psi'_t(x)] \nabla g_t(\psi_t(x))$, where $[\psi'_t]$ denotes the $N \times N$ Jacobian matrix  of $\psi_t$  with respect to $x$. Thus,  $\nabla g_t(\psi_t(x)) =  [ \psi'_t(x) ]^{-1} \nabla u_0(x)$, and consequently, we have $|\nabla g_t(\psi_t(x))|_A^p = | [ \psi'_t(x) ]^{-1} \nabla u_0(x)|^p_A $. Then, under the transformation $y=\psi_t(x)$, we have
\begin{align}\label{deriv-4}
      \int_{\Omega} |\nabla g_t(y)|^p_A \d y &=  \int_{\Omega}  | [ \psi'_t(x) ]^{-1} \nabla u_0(x)|^p_A |J\psi_t| \d x \nonumber\\
       &= \int_{\Omega}  | [ \psi'_t(x) ]^{-1} \nabla u_0(x)|^p_A ( 1 + t \operatorname{div} W + o(t^2)) \d x \nonumber\\
       & =\int_{\Omega} \left|\left(I-t W' + o(t^2) \right) \nabla u_0 \right|_A^p \left( 1 + t \operatorname{div} W + o(t^2)\right) \d x \nonumber\\
      % & =\int_{\Omega} | \nabla u_0^T -  t \nabla u_0^T W' + o(t^2) |_A^p ( 1 + t \operatorname{div} W + o(t^2)) \d x \\
       &= \int_{\Omega} \left[|\nabla u_0|^p_A - t p |\nabla u_0|^{p-2}_A \left< A \nabla u_0 , W' \nabla u_0 \right> \right] (1 + t \operatorname{div} W ) \d x + o(t^2)\nonumber\\
       &= \int_{\Omega}|\nabla u_0|_A^p \d x + t \int_{\Omega}|\nabla u_0|_A^p  \operatorname{div} W \d x \nonumber\\
       &\quad - tp\int_{\Omega} |\nabla u_0|_A^{p-2} \left< A\nabla u_0, W' \nabla u_0 \right> \d x + o(t^2).
\end{align}
Here in the above  the first order Taylor expansion $|X+tY|^p_A = |X|^p_A + tp |X|^{p-2}_A \left< AX, Y \right> + o(t^2)$ has been used for vectors $X = \nabla u_0$ and $Y  = - W' \nabla u_0$. Then, using \eqref{deriv-2}, \eqref{deriv-3} and \eqref{deriv-4} in \eqref{deriv-1}, we deduce that
\begin{align*}
    \lambda_1(V_t,m_t) &( 1  + t \int_{\Omega} m |u_0|^p\operatorname{div} W\d x )\nonumber\\
    &\leq \int_{\Omega}|\nabla u_0|^p_A \d x+ \int_{\Omega} V |u_0|^p \d x +  t \int_{\Omega} (|\nabla u_0|^p_A  + V|u_0|^p )\operatorname{div} W\d x \nonumber\\
    & \hspace{2.5cm}- tp \int_{\Omega} |\nabla u_0|^{p-2}_A \left< A \nabla u_0, W' \nabla u_0 \right> \d x + o(t^2) \nonumber  \\
    & = \lambda_1(V,m) + t\int_{\Omega} \left[ |\nabla u_0|_A^p + V|u_0|^p\right]\operatorname{div} W d x \nonumber\\
    &\hspace{2.5cm}- tp \int_{\Omega} |\nabla u_0|^{p-2}_A \left< A \nabla u_0, W' \nabla u_0 \right> \d x + o(t^2).
\end{align*}
Then, it follows that
    \begin{align}\label{diff-3}
    &\frac{\lambda_1(V_t,m_t) - \lambda_1(V,m)}{t} \nonumber\\
     &\hspace{1cm}\leq \int_{\Omega}\left[|\nabla u_0|^p_A  + V|u_0|^p \right]\operatorname{div} W\d x - p \int_{\Omega} |\nabla u_0|_A^{p-2} \left< A \nabla u_0, W' \nabla u_0 \right> \d x  \nonumber \\
    &\hspace{6.5cm}-\lambda_1(V_t,m_t) \int_{\Omega} m |u_0|^p\operatorname{div} W\d x + o(t). 
\end{align}
    Similarly, using  $f_t := u_t \circ \psi_t$ as an admissible function in the definition of $\lambda_1(V,m)$, one can find that 
\begin{align}\label{diff-4}
    &\frac{\lambda_1(V_t,m_t) - \lambda_1(V,m)}{t} \nonumber\\
    &\hspace{1cm}\geq \int_{\Omega} \left[|\nabla f_t|^p_A  + V|f_t|^p \right]\operatorname{div} W\d x  - p \int_{\Omega} |\nabla f_t|_A^{p-2} \left< A \nabla f_t, W' \nabla f_t\right> \d x \nonumber  \\
    & \hspace{6.5cm} -\lambda_1(V,m) \int_{\Omega} m |f_t|^p\operatorname{div} W\d x + o(t). 
\end{align}
Since the vector field $W$ is smooth and compactly supported, hence $\operatorname{div} W \in L^\infty(\Omega)$. Also, by the Lemma \ref{convergence-u}, $f_t \to u_0$ in $\Wkp$ as $t \to 0^+$. Hence,
 \begin{align}\label{diff-8}
     \int_{\Omega}| \nabla f_t |_A^p \operatorname{div} W \d x \to \int_{\Omega}| \nabla u_0 |_A^p \operatorname{div} W \d x \quad \text{ as } t \to 0^+.
 \end{align}
 Further, applying the Lemma \ref{convergence-3}, we have
 \begin{align}\label{diff-6}
     \int_{\Omega} V|f_t|^p \operatorname{div} W \d x \to \int_{\Omega} V|u_0|^p \operatorname{div} W \d x,
 \end{align}
  and 
  \begin{align}\label{diff-7}
      \int_{\Omega} m|f_t|^p \operatorname{div} W \d x \to \int_{\Omega} m|u_0|^p \operatorname{div} W \d x 
  \end{align}
 as $ t \to 0^+$.  
Now, consider a map $G : \Wkp \to L^1(\Omega)$ defined by $$G(u) := |\nabla u|^{p-2}_A \left< A \nabla u, W'\nabla u \right>.$$ Then, $ |G(u)|  \leq \Theta^{-p}_\Omega \|W'\|_\infty  |\nabla u|^p$.
Suppose that $u_n \to u$ in $\Wkp$. Then $|\nabla u_n|^p \to |\nabla u|^p$ in $L^1(\Omega)$, and up to a subsequence (still denoted by $\{\nabla u_n\}$),  $\nabla u_n \to \nabla u$ a.e. in $\Omega$. This implies that $G(u_n) \to G(u)$ a.e. in $\Omega$. 
Hence by the  Dominated convergence theorem, $\|G(u_n)-G(u)\|_1 \to 0 $ as $n \to \infty$. This shows that $G$ is continuous in $\Wkp$. Consequently, we have $G(f_t) \to G(u_0)$ in $L^1(\Omega)$ as $t \to 0^+$ i.e., 
 \begin{align}\label{diff-9}
     \int_{\Omega} |\nabla f_t|_A^{p-2} \left< A \nabla f_t, W' \nabla f_t\right> \d x \to  \int_{\Omega} |\nabla u_0|_A^{p-2} \left< A \nabla u_0, W' \nabla u_0 \right> \d x \text{ as } t \to 0^+.
 \end{align}
 Finally, letting $ t \to 0^+$ in \eqref{diff-3} and  \eqref{diff-4}, and applying  \eqref{diff-8}, \eqref{diff-6}, \eqref{diff-7}, together with \eqref{diff-9}, we deduce that
\begin{align*}
   \frac{\d \lambda_1(t) }{\d t} \Bigg|_{t=0} &=   \lim_{t \to 0^+}  \frac{\lambda_1(V_t,m_t) - \lambda_1(V,m)}{t}\\
   & =\int_{\Omega} \left[|\nabla u_0|^p_A  + V|u_0|^p \right]\operatorname{div} W\d x  - \lambda_1(V,m) \int_{\Omega}m |u_0|^p  \operatorname{div} W \d x\\
    &\hspace{3cm}- p \int_{\Omega}|\nabla u_0|^{p-2}_A \left<A \nabla u_0, W' \nabla u_0 \right> \d x.
\end{align*}
This concludes the result.   
\end{proof}

%\noindent \textbf{Conflict of interest statement:} On behalf of the authors, the corresponding author states that there is no conflict of interest.
%\newline
%\textbf{Data availability statement:} Data sharing does not apply to this article as no datasets were generated or analyzed during the current study.
	\section*{Acknowledgement}
The authors express gratitude to Prof. Yehuda Pinchover for his helpful suggestions concerning the Morrey spaces. The author Pradeep Gupta acknowledges the financial support received from the DST, Govt. of India, through the DST Inspire Fellowship (Grant No. DST/INSPIRE Fellowship/2024/IF240239). R.Kr. Giri thanks the DST-FIST, Govt. of India for providing financial support under the scheme “Fund for Improvement of Science and Technology” (FIST - No. SR/FST/MS-I/2018/24). S. Ghosh gratefully acknowledges financial support for this research under the ARG-MATRICS grant No. ANRF/ARGM/2025/001570/MTR, ANRF, Govt. of India.

% \cite{Arias02}
% \cite{Henrot2019,Henrot2006,Veronica2025,Noris2026}
% \cite{drabek97}
% \cite{Glowinski75}\cite{dold1975}\cite{lakshmi2026}\cite{Glowinski03,gilboa18,Bogosel18,Elmoataz15,Szulkin88,Otani88,Strue08,keshavan12,Cuesta09,Giri20,Giri23,Pinchover16,Otani88,Cuesta01,T.Anoop11,H.Y.X95,W.Alle95,W.Alle98,K.J.brn90,P.Lind90,Garc87,Brezis11,Pucci07,Pel76}
% 	\cite{Ziemer97}\cite{Anane88} \cite{pinchover2008}
%     \cite{Anle06}

	\bibliographystyle{plain}
	\bibliography{Reference}
\end{document}